\documentclass[11pt]{article}
\usepackage[utf8]{inputenc}
\usepackage[english]{babel}
\usepackage{graphicx}
\usepackage{amsmath}
\usepackage{amssymb}
\usepackage{amsfonts}
\usepackage{amsthm}
\usepackage[left=2.1cm,top=1.5cm,right=2.1cm, bottom=2.2cm,letterpaper]{geometry}
\usepackage{mathrsfs}
\usepackage{caption}
\usepackage{subcaption}
\usepackage{latexsym}
\usepackage{xfrac}
\usepackage{tcolorbox}
\usepackage{enumitem}
\usepackage{float}
\usepackage{hyperref}
 \usepackage{url}
\usepackage[toc]{appendix}
\newtheorem{theorem}{Theorem}[section]
\newtheorem{corollary}{Corollary}[section]

\newtheorem{definition}{Definition}[section]
\newtheorem{proposition}{Proposition}[section]

\newtheorem{remark}{Remark}[section]

\numberwithin{equation}{section}
\numberwithin{figure}{section}

\newcommand{\N}{\mathbb{N}}
\newcommand{\R}{\mathbb{R}}

\newcommand{\Rey}{\rm Re}

\makeatletter
\renewcommand*\env@matrix[1][*\c@MaxMatrixCols c]{%
	\hskip -\arraycolsep
	\let\@ifnextchar\new@ifnextchar
	\array{#1}}
\makeatother

\def\eq#1{(\ref{#1})}
\def\neweq#1{\begin{equation}\label{#1}}
\def\endeq{\end{equation}}

\begin{document}

\title{A connection between symmetry breaking for Sobolev minimizers \\ and stationary Navier-Stokes flows past a circular obstacle}

\author{Filippo GAZZOLA -- Gianmarco SPERONE -- Tobias WETH}
\date{}
\maketitle
\vspace*{-6mm}

\begin{abstract}
\noindent
Fluid flows around a symmetric obstacle generate vortices which may lead to symmetry breaking of the streamlines. We study this phenomenon for planar
viscous flows governed by the stationary Navier-Stokes equations with constant inhomogeneous Dirichlet boundary data in a rectangular channel containing a circular obstacle.
In such (symmetric) framework, symmetry breaking is strictly related to the appearance of multiple solutions.
Symmetry breaking properties of some Sobolev minimizers are studied and explicit bounds on the boundary velocity (in terms of the length and
height of the channel) ensuring uniqueness are obtained after estimating some Sobolev embedding constants and constructing a suitable solenoidal extension of the boundary data.
We show that, regardless of the solenoidal extension employed, such bounds converge to zero at an optimal rate as the length of the channel tends to infinity.\par\noindent
\end{abstract}


\section{Introduction}

We are interested in showing a connection between some purely theoretical tools in functional analysis, such as Sobolev minimizers, and some problems in fluid
mechanics, such as the appearance of vortices in flows around an obstacle. 
It is well-known that, if a quickly moving fluid hits an obstacle, a vortex shedding appears behind
the body. This happens also for flows which are laminar before hitting the obstacle.

\begin{figure}[H]
\begin{center}
\includegraphics[width=168mm]{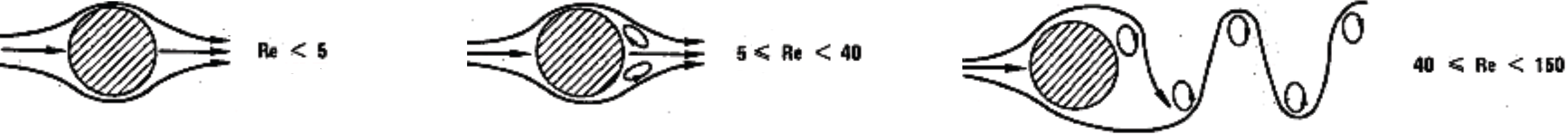}\end{center}
\vspace*{-5mm}
\caption{Regimes of flow across a circular cylinder, taken from \cite{Lienhard}.}\label{vortexshedding}
\end{figure}

In literature, the most documented fluid-structure interaction experiment is that of a circular pendulum (cylinder) in a water flow, where
the interaction consists both in the generation of vortices (structure modifying the fluid) and in the appearance of lift forces
(exerted by the fluid on the structure), namely forces transversal to the flow; we refer to
\cite{CoutanceauD,Lienhard,Morkovin64,Schlichting,Tritton59,Williamson}
and to \cite{denis} for more references. When a cylinder of diameter $d$ is placed with its axis normal
to the flow having an upstream constant speed $U$, simple experiments enable to highlight turbulence, vortex shedding, and the appearance
of the lift force, see e.g.\ \cite{wolfgang4}. If the cylinder is long compared with $d$, this experiment can be modeled with a disk in a 2D flow.
One can then vary the value of the speed $U$ and observe that the flow pattern depends on the Reynolds number ($\text{Re}= U d / \eta$,
with $\eta$ being the viscosity of the fluid).
The pictures in Figure \ref{vortexshedding} summarize the observations and the symmetry breaking according to the Reynolds number.
In the left picture ($\Rey<5$) the flow has a double symmetry, both vertical and horizontal, and the regime is that of {\em unseparated flow}.
In the middle picture ($5<\Rey<40$) there is symmetry breaking in the horizontal direction, with the appearance of a pair of F\"oppl vortices
in the wake. In the right picture ($40<\Rey<150$) there is symmetry breaking in the vertical direction leading to a {\em laminar vortex street}.
For larger $\Rey$ the picture becomes more and more disordered with the appearance of {\em turbulence and  chaos}, see \cite{Lienhard}. In the present paper we will see that a similar phenomenon occurs for the minimizers of the Sobolev ratio as the fluid domain tends to enlarge. \par

\begin{figure}[H]
	\begin{center}
		\includegraphics[height=45mm,width=99mm]{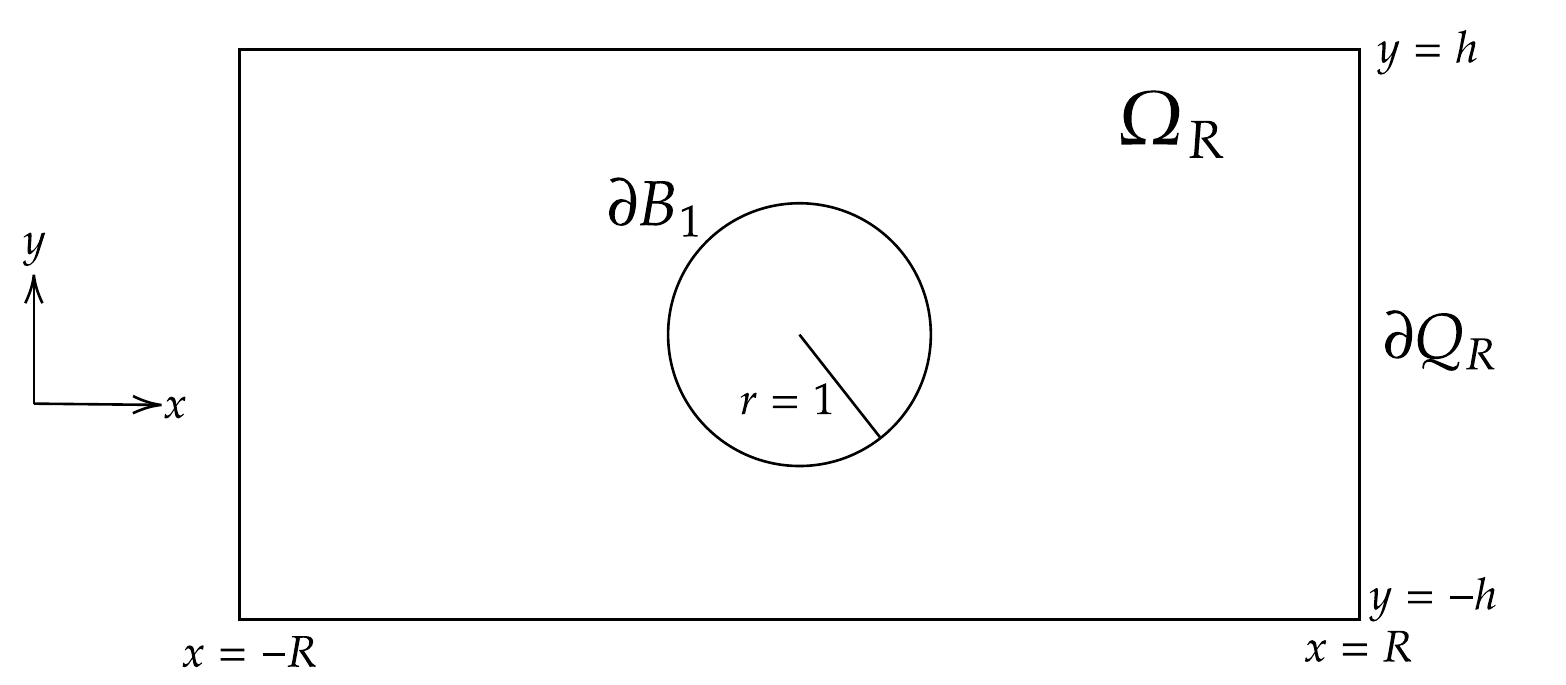}
	\end{center}
	\vspace*{-4mm}
	\caption{The pierced rectangle $\Omega_{R}$ in \eqref{rectangle}}\label{dom1}
\end{figure}

Through the Navier-Stokes equations we model the steady motion of a viscous incompressible fluid in a large horizontal channel containing the unit disk $B_1$ (circular
obstacle of diameter $d=2$) under constant inflow, no-slip boundary conditions and no external force, see \cite{amick1977steady,bonheure2020equilibrium,morimoto2004remark}. To this end we take a pierced rectangle with
variable length $2R$ in the $x$-direction and (fixed) given height $2h$ in the $y$-direction:
\begin{equation} \label{rectangle}
Q_{R} \doteq (-R,R) \times (-h,h) \qquad \text{and} \qquad \Omega_{R} \doteq Q_{R }\setminus \overline{B_{1}}\, ,
\end{equation}
with $R>h>1$, see Figure \ref{dom1}.

We study the following boundary-value problem related to the stationary Navier-Stokes equations:
\begin{equation}\label{nsstokes0}
\left\{
\begin{aligned}
& -\eta\Delta u+(u\cdot\nabla)u+\nabla p=0,\ \quad  \nabla\cdot u=0 \ \ \mbox{ in } \ \ \Omega_{R}, \\[3pt]
& u = (U,0) \ \text{ on } \ \partial Q_R, \qquad u = (0,0) \ \text{ on } \ \partial B_{1} \, .
\end{aligned}
\right.
\end{equation}
In \eqref{nsstokes0}, $\eta>0$ is the (constant) kinematic viscosity of the fluid, $u : \Omega_{R} \longrightarrow \mathbb{R}^2$ is the velocity vector field and $p : \Omega_{R} \longrightarrow \mathbb{R}$ is the scalar pressure. The boundary conditions in \eqref{nsstokes0}$_2$ model the inflow/outflow of the fluid across the outer boundary $\partial Q_R$ with a given constant speed $U > 0$ (in the $x$-direction), and with no-slip condition on the obstacle $B_1$ where viscosity yields zero velocity of the flow. Notice that the compatibility condition
\begin{equation} \label{zeroflux2}
\int_{\Omega_R}\nabla \cdot u =\int_{\partial\Omega_R}u \cdot \nu = \int_{\partial Q_R} (U,0) \cdot \nu = 0
\end{equation}
is satisfied for the solutions of \eq{nsstokes0}. In  \eqref{zeroflux2}, $\nu \in \mathbb{R}^2$ denotes the outward unit normal to $\Omega_R$.\par
In \cite{gazspe} (see also \cite{kochuniqueness,denis,gazspefra,filippoclara,sperone2021steady}) it is shown that, in a symmetric framework,
the most important threshold for Re is the one yielding symmetry breaking for solutions of \eq{nsstokes0}, since it corresponds to the appearance of the lift force exerted from the fluid on the obstacle. But symmetry breaking also means multiplicity of the equilibrium configurations of the fluid. Sufficient conditions for uniqueness of solutions for \eq{nsstokes0} are usually obtained through smallness assumptions on the inflow and they involve purely analytical quantities such as Sobolev embedding constants. This raises the natural question
\begin{center}
In a symmetric framework, is symmetry breaking for solutions of \eq{nsstokes0} (that is, nonuniqueness) \\
related to symmetry breaking of some minimizers for Sobolev ratios?
\end{center}
In the present paper, we start tackling this problem and we show that, indeed, there is a connection between symmetry breaking in the stationary Navier-Stokes equations and Sobolev minimizers. This requires the determination of explicit bounds both for some Sobolev embedding constants and for unique solvability of the Navier-Stokes equations. With these bounds at hand we exhibit a qualitative connection between the symmetry breaking in Figure \ref{vortexshedding} as $\text{Re}$ increases, and in Sobolev minimizers within the domain $\Omega_{R}$ as $R \to \infty$. This connection is properly explained in Remarks \ref{conn1} and \ref{conn2} below.

\section{Bounds for some Sobolev embedding constants}\label{bondsSob}

We consider the Sobolev space $H_{0}^{1}(\Omega_{R})$ (for both scalar functions and vector fields) and the space of scalar  $H^{1}(\Omega_{R})$-functions that are constant on $\partial Q_{R}$ and vanish on $\partial B_{1}$, which is a proper connected part of
$\partial\Omega_{R}$:
$$
H^1_c(\Omega_{R}) \doteq \{v\in H^1(\Omega_{R}) \ | \ v=0 \ \ \mbox{on} \ \ \partial B_{1} \, , \ v \ \text{is constant on} \ \partial Q_{R} \}\, .
$$
Since $|\partial B_{1}|_{1} = 2 \pi >0$ (the 1D-Hausdorff measure),
the Poincaré inequality holds in $H^1_c(\Omega_{R})$, which means that $v\mapsto\|\nabla v\|_{L^2(\Omega_{R})}$ is a norm on $H^1_c(\Omega_{R})$. This space may be rigorously characterized by using the relative capacity potential $\psi \in H^{1}_{0}(\Omega_{R})$ of $B_{1}$ with respect to $Q_R$, defined as the function $\psi \in H^{1}(\Omega_{R})$ such that
$$
\begin{aligned}
& \Delta\psi=0 \quad \text{in} \ \ \Omega_{R}, \qquad \psi = 0 \quad \text{on} \ \ \partial Q_{R},\qquad \psi=1 \quad \text{in} \ \ B_{1} \, , \\[6pt]
& \hspace{0.5cm} \|\nabla\psi\|^2_{L^2(\Omega_{R})} = \min_{v\in H_{0}^{1}(Q_{R})} \left\{  \int_{\Omega_{R}} |\nabla v|^2 \ \Big| \ v=1 \ \text{ in } \ \overline{B_{1}} \, \right\} \, .
\end{aligned}
$$
Therefore, $H^1_c(\Omega_{R})$ has the characterization
$$
H^1_c(\Omega_{R})=H^1_0(\Omega_{R})\oplus \R(\psi-1) \quad \text{with} \quad H^1_0(\Omega_{R})\perp\R(\psi-1)\, ,
$$
so that $H^1_0(\Omega_{R})$ has codimension 1 within $H^1_c(\Omega_{R})$ and the ``missing dimension'' is spanned by the function $\psi-1$. To see this,
determine the orthogonal complement of $H^1_0(\Omega_{R})$ within $H^1_c(\Omega_{R})$ as follows:
\begin{equation} \label{ortcom}
\begin{aligned}
v\in H^1_0(\Omega_{R})^\perp\qquad & \Longleftrightarrow \quad v\in H^1_c(\Omega_{R}) \quad \text{and} \quad (\nabla v, \nabla w)_{L^{2}(\Omega_{R})} =0 \quad \forall w\in H^1_0(\Omega_{R}) \\[6pt]
& \Longleftrightarrow\quad v\in H^1_c(\Omega_{R}) \quad \text{and} \quad \langle\Delta v,w\rangle_{\Omega_{R}}=0 \quad \forall w\in H^1_0(\Omega_{R}) \, ,
\end{aligned}
\end{equation}
so that $v$ is weakly harmonic and, since $v\in H^1_c(\Omega_{R})$, it is necessarily a real multiple of $\psi-1$. In \eqref{ortcom}, $\langle \cdot , \cdot \rangle_{\Omega_{R}}$ denotes the duality product between $H^{-1}(\Omega_{R})$ and $H^1_0(\Omega_{R})$.
\par
We then define the Sobolev constant of the embedding $H_0^{1}(\Omega_{R}) \subset L^{4}(\Omega_{R})$ as
\begin{equation} \label{sobolevconstants1}
\mathcal{S}_R = \min_{v \in H_0^{1}(\Omega_{R}) \setminus \{0\}} \ \frac{\|\nabla v\|^2_{L^2(\Omega_{R})}}{\|v\|^{2}_{L^{4}(\Omega_{R})}} \, ,
\end{equation}
so that $\mathcal{S}_{R} > 0$ and
\begin{equation} \label{sobolevconstants11}
\mathcal{S}_{R} \, \|v\|^{2}_{L^{4}(\Omega_{R})} \leq \|\nabla v\|^2_{L^2(\Omega_{R})} \qquad \forall v \in H_0^{1}(\Omega_{R}) \, .
\end{equation}

The main goal of this section is to provide explicit lower and upper bounds for $\mathcal{S}_{R}$, in terms of $R$ and $h$. For later use we introduce
\neweq{besselmu}
\mu_0=\mbox{the first zero of the Bessel function of first kind of order zero}\approx2.40483\, .
\endeq

\begin{theorem}\label{boundsobolevrec}
For $\mathcal{S}_R$ as in \eqref{sobolevconstants1} we have
\begin{equation} \label{lowerboundrec}
\pi \sqrt{ \dfrac{3}{2} } \, \max \left\{ \dfrac{\sqrt{\pi}}{2} \dfrac{\sqrt{R^2 + h^2}}{R h} ,\frac{\mu_0}{\sqrt{4Rh - \pi}}\right\} \leq \mathcal{S}_{R} \leq \dfrac{1}{2} \sqrt{\dfrac{\pi}{2 \kappa}} (2 h^2 \log(h)^2 - 2 h^2 \log(h) + h^2 - 1) \, ,
\end{equation}
where
$$
\begin{aligned}
\kappa & \doteq - \dfrac{1}{324} + \dfrac{96 h}{3125} - \dfrac{9 h^2}{64} + \dfrac{32 h^3}{81} - \dfrac{3 h^4}{4} + \dfrac{7580461 h^6}{16200000} \\[5pt]
& \hspace{4.5mm} - \dfrac{66801 h^6 \log(h) - 46690 h^6 \log(h)^2 + 17400 h^6 \log(h)^3 -
	3000 h^6 \log(h)^4}{90000} \, .
\end{aligned}
$$
In particular, for any scalar or vector function $w\in H^1_0(\Omega_{R})$ one has
	\neweq{ineqL4H1u2rec}
	\|w\|_{L^4(\Omega_{R})}^2 \leq \dfrac{1}{\pi} \sqrt{ \dfrac{2}{3} } \, \min\left\{ \dfrac{2}{\sqrt{\pi}} \dfrac{R h}{\sqrt{R^2 + h^2}} ,\frac{\sqrt{4Rh - \pi}}{\mu_0}\right\}\, \, \|\nabla w\|_{L^2(\Omega_{R})}^2\, .
	\endeq

\end{theorem}
\noindent
\begin{proof}
 For a scalar function $w\in H^1_0(\Omega_{R})$, after combining the Gagliardo-Nirenberg inequality in $\R^2$ given by del Pino-Dolbeault \cite[Theorem 1]{delpino} with some H\"older inequality, the following inequality was obtained in \cite[Theorem 2.3]{gazspe}:
\neweq{holder3}
\|w\|_{L^4(\Omega_{R})}^2\le \sqrt{\frac{2}{3\pi}} \, \|\nabla w\|_{L^2(\Omega_{R})}\|w\|_{L^2(\Omega_{R})}\qquad\forall w\in H^1_0(\Omega_{R})\, .
\endeq
In fact, \eqref{holder3} improves previous results by Ladyzhenskaya \cite{Lady59} (see also \cite[Lemma 1, p.8]{ladyzhenskaya1969mathematical}) and Galdi \cite[(II.3.9)]{galdi2011introduction}.	
\par
We then notice that
	$$
	-\Delta \left( \cos \left( \frac{\pi x}{2R} \right) \cos \left( \frac{\pi y}{2h} \right) \right) = \dfrac{\pi^{2}}{4} \dfrac{R^2 + h^2}{R^2 h^2} \cos \left( \frac{\pi x}{2R} \right) \cos \left( \frac{\pi y}{2h} \right) \qquad \forall (x,y) \in Q_{R} \, ;
	$$
	such eigenfunction for the Dirichlet-Laplace operator is positive in $Q_{R}$ and vanishes on $\partial Q_R$. Therefore, by applying the Poincar\'e inequality in $Q_{R}$ to $w \in H^1_0(Q_{R})$, we obtain
	\begin{equation} \label{poincrec1}
	\|w\|_{L^2(\Omega_{R})} \leq \dfrac{2}{\pi} \dfrac{R h}{\sqrt{R^2 + h^2}} \, \|\nabla w\|_{L^2(\Omega_{R})} \, .
	\end{equation}
	
	In order to finish the proof for scalar functions we apply a symmetrization argument. Let $\Omega^* \subset \mathbb{R}^{2}$ be a disk having the same measure as $\Omega_{R}$, and thus its radius is
	$$
	R_{*}=\sqrt{\dfrac{4Rh - \pi}{\pi}} \, .
	$$
	Since the Poincaré constant of $B_1$ is given by $\mu_{0}^{2}$ (see \eqref{besselmu}), by rescaling, the Poincaré constant of $\Omega^*$ is $\mu_{0}^2/R_{*}^2$. In view of the
	Faber-Krahn inequality \cite{faber1923beweis,krahn1925rayleigh} this implies that
	$$
	\min_{v\in H^1_0(\Omega_{R})}\ \frac{\|\nabla v\|_{L^2(\Omega_{R})}}{\|v\|_{L^2(\Omega_{R})}}\, \ge\, \min_{v\in H^1_0(\Omega^*)}\ \frac{\|\nabla v\|_{L^{2}(\Omega^{*})}}{\|v\|_{L^{2}(\Omega^{*})}}
	=\frac{\mu_{0}}{R_{*}}\, ,
	$$
	and therefore
	\begin{equation} \label{poincrec2}
	\|w\|_{L^2(\Omega_{R})} \leq \frac{R_{*}}{\mu_{0}} \|\nabla w \|_{L^2(\Omega_{R})} = \frac{1}{\mu_{0}} \sqrt{\dfrac{4Rh - \pi}{\pi}} \, \|\nabla w \|_{L^2(\Omega_{R})} \, .
	\end{equation}
	The lower bound in \eqref{lowerboundrec} is reached after inserting \eqref{poincrec1} and \eqref{poincrec2} into \eqref{holder3}, which then yields \eqref{ineqL4H1u2rec}. \par The upper bound in \eqref{lowerboundrec} is obtained by testing the quotient \eqref{sobolevconstants1} with
	\begin{equation} \label{testrec}
	X_{0}(x,y) \doteq \dfrac{1}{2} \left( h-\sqrt{x^2 + y^2} \right) \log(x^2 + y^2) \qquad \forall (x,y) \in  \Omega_{R} \, ,
	\end{equation}
	which, if extended by zero for $x^2 + y^2 \geq h^2$, becomes an element of $H^1_0(\Omega_{R})$.

\begin{figure}[H]
	\begin{center}
		\includegraphics[height=65mm,width=100mm]{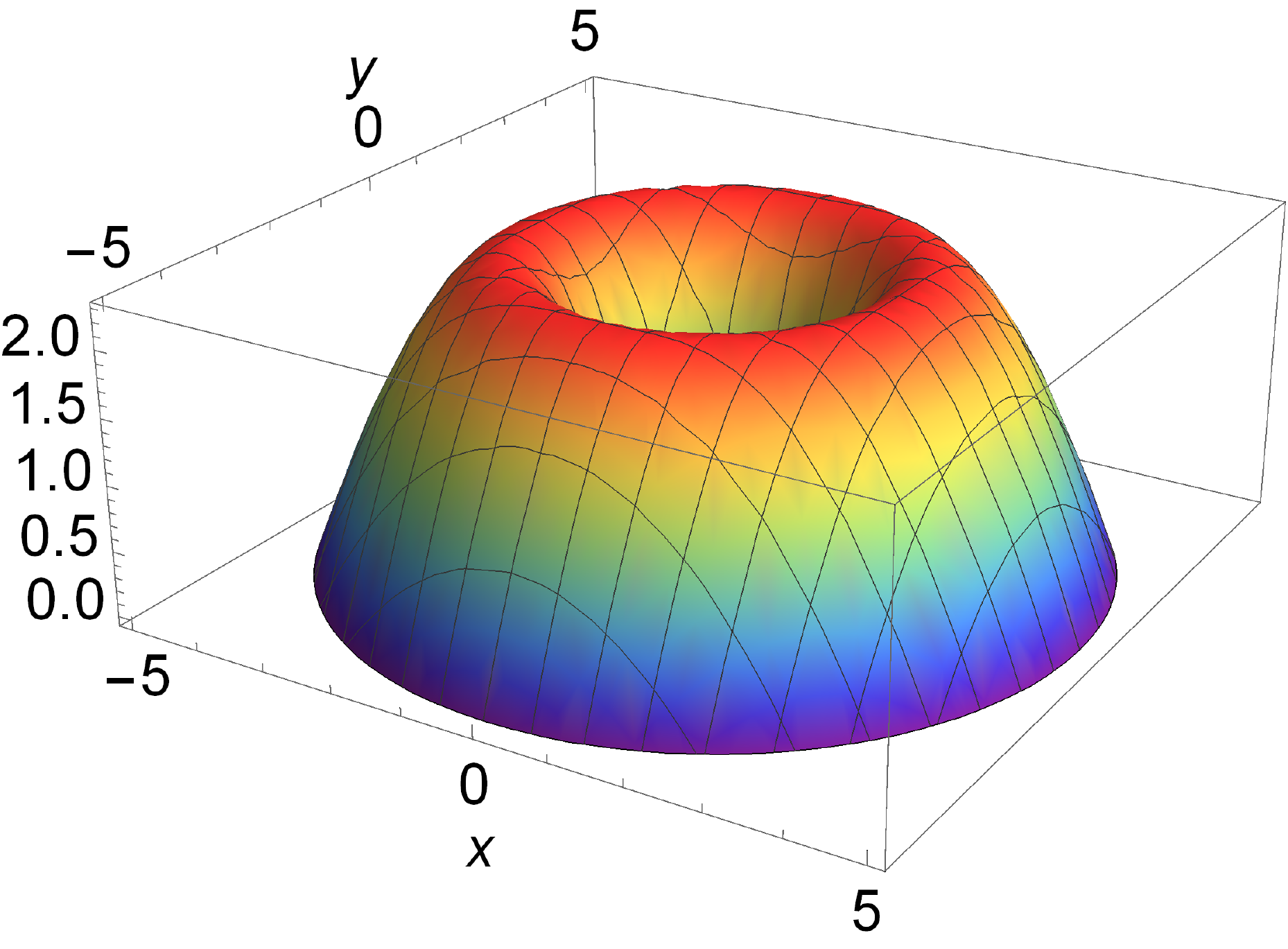}
	\end{center}
	\vspace*{-5mm}
	\caption{Graph of the function $X_{0}$ in \eqref{testrec} for $h=5$.}\label{sobrecub}
\end{figure}
	
To conclude we notice that inequality \eqref{sobolevconstants11} is also valid for vector functions (with the same constant): if
$v = (v_1, v_2) \in H_{0}^{1}(\Omega_{R})$ is a vector field, by the Minkowski inequality we get
$$
\begin{aligned}
\|v\|^{4}_{L^4(\Omega_{R})} & = \left\| \, | v_{1} |^{2} + | v_{2} |^{2} \, \right\|^{2}_{L^{2}(\Omega_{R})} \leq \left( \|v_{1}\|^{2}_{L^4(\Omega_{R})} + \|v_{2}\|^{2}_{L^4(\Omega_{R})} \right)^{2} \\[3pt]
& \leq \mathcal{S}_{R}^{-2} \left( \| \nabla v_{1}\|^{2}_{L^2(\Omega_{R})} + \| \nabla v_{2}\|^{2}_{L^2(\Omega_{R})}  \right)^{2} =
\mathcal{S}_{R}^{-2} \| \nabla v \|^{4}_{L^2(\Omega_{R})} \, .
\end{aligned}
$$
Therefore, \eq{sobolevconstants11} also holds for vector functions.	
\end{proof}

\begin{remark} \label{sobquad}
A careful look at the proof of Theorem \ref{boundsobolevrec} shows that \eqref{holder3} and \eqref{poincrec1} also hold for functions in $H^{1}_{0}(Q_R)$, namely
\begin{equation} \label{sobquad1}
\|w\|_{L^4(\Omega_{R})}^2 \leq \dfrac{2}{\pi} \sqrt{ \dfrac{2}{3 \pi} } \dfrac{R h}{\sqrt{R^2 + h^2}} \|\nabla w\|_{L^2(\Omega_{R})}^2 \qquad \forall w \in H^{1}_{0}(Q_R).
\end{equation}
\end{remark}

\begin{remark} \label{remrec}
The upper bound in \eqref{lowerboundrec} is obtained by using the function $X_0$ in \eqref{testrec}. For $R < 2h + 1$, a smaller upper bound may be derived through the Sobolev constant of the rectangle $(1,R) \times (-h,h)$, which contains the disc $D_{R, h}$ of radius $h$ centered at $\left( \frac{R+1}{2}, 0 \right)$: then one considers the function
$$
X_{1}(x,y) \doteq  h^{2} - \left( x - \dfrac{R+1}{2} \right)^{2} - y^{2} \qquad \forall (x,y) \in  D_{R, h} \, ,
$$
which gives the estimates
\begin{equation} \label{soblimit0}
\pi \sqrt{ \dfrac{3}{2} } \, \max \left\{ \dfrac{\sqrt{\pi}}{2} \dfrac{\sqrt{R^2 + h^2}}{R h} ,\frac{\mu_0}{\sqrt{4Rh - \pi}}\right\} \leq \mathcal{S}_{R} \leq \dfrac{2 \sqrt{5 \pi}}{h} \qquad \text{for} \ R < 2h+1 \, .
\end{equation}
Several further trials were executed with different functions in $H^1_0(\Omega_{R})$, but \eqref{testrec} appears to give a good upper bound for $\mathcal{S}_{R}$ for a wide range of values of $h >1$. From \eqref{lowerboundrec} we observe that
\begin{equation} \label{soblimit1}
\dfrac{\pi}{2} \sqrt{ \dfrac{3 \pi}{2} } \dfrac{1}{h} \leq \lim\limits_{R \to \infty} \mathcal{S}_{R} \leq \dfrac{1}{2} \sqrt{\dfrac{\pi}{2 \kappa}} (2 h^2 \log(h)^2 - 2 h^2 \log(h) + h^2 - 1) \qquad \forall h >1 \, .
\end{equation}
The plot in Figure \ref{logwins} below compares the lower and upper bounds in \eqref{soblimit0}-\eqref{soblimit1} as functions of $h >1$, as $R \to \infty$. A simple computation shows that the ratio between the upper and lower bounds in \eqref{soblimit1} tends to $2 \sqrt{10}/ \pi \approx 2.01317$ as $h \to \infty$.
\begin{figure}[H]
	\begin{center}
		\includegraphics[height=70mm,width=160mm]{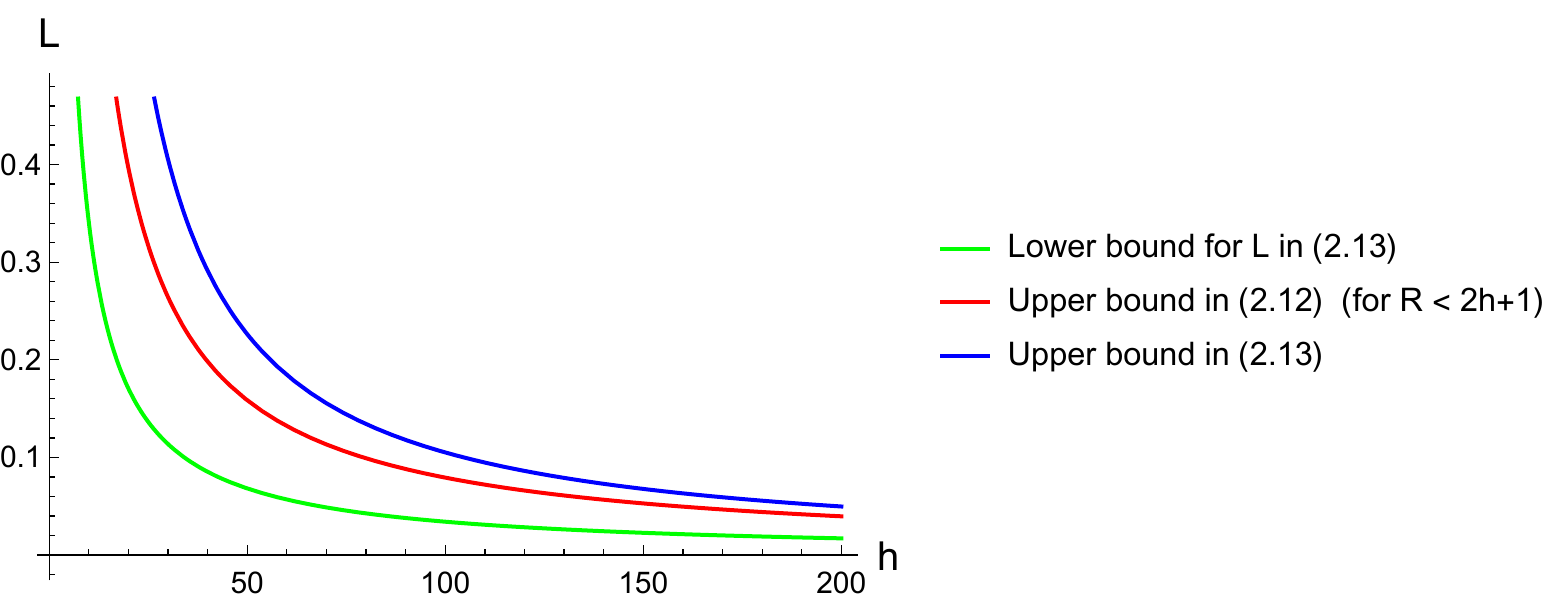}
	\end{center}
	\vspace*{-6mm}
	\caption{Behavior of the lower and upper bounds for $ L \doteq \lim\limits_{R \to \infty} \mathcal{S}_{R} $ as functions of $h > 1$.}
	\label{logwins}
\end{figure}
\noindent
Later, in Proposition \ref{betterbound} we shall improve the upper bound in \eqref{soblimit1} after verifying that $\lim\limits_{R \to \infty} \mathcal{S}_{R}$ coincides with the best Sobolev constant of the infinite strip $\mathbb{\R} \times (-h,h)$, for every $h >1$ (see Theorem \ref{tobias}).
\end{remark}

\section{Symmetry breaking of some Sobolev minimizers}\label{symbreakSob}


In this section we consider the subspace $H_{0,e}^{1}(\Omega_{R})$ of functions $u \in H_0^{1}(\Omega_{R})$ which are even in $x$, that is, $u(x,y)=u(-x,y)$ for every $(x,y) \in \Omega_R$, and we set
\begin{equation} \label{sobolevconstants1-R-even}
\mathcal{S}_{R,e} = \min_{v \in H_{0,e}^{1}(\Omega_{R}) \setminus \{0\}} \ \frac{\|\nabla v\|^2_{L^2(\Omega_{R})}}{\|v\|^{2}_{L^{4}(\Omega_{R})}} \, .
\end{equation}

We prove here the following symmetry breaking result.

\begin{theorem}\label{symm-breaking} \label{tobias}
  There exists $R_0>0$ with the property that $\mathcal{S}_{R,e} > \mathcal{S}_{R}$ for $R>R_0$. Hence every minimizer of the Sobolev quotient \eqref{sobolevconstants1} is a non-even function in the $x$-variable for $R>R_0$.
\end{theorem}
\noindent
\begin{proof} Consider the infinite strip $\Omega_\infty\doteq \R \times (-h,h)$ and put
\begin{equation} \label{sobolevconstants1-star}
\mathcal{S}_\infty = \inf_{v \in H_0^{1}(\Omega_\infty) \setminus \{0\}} \ \frac{\|\nabla v\|^2_{L^2(\Omega_\infty)}}{\|v\|^{2}_{L^{4}(\Omega_\infty)}}. 
\end{equation}
Since for $R'>R>1$ we have the inclusions $H_0^{1}(\Omega_R) \subset H_0^1(\Omega_{R'}) \subset H_0^{1}(\Omega_\infty)$ induced by trivial extension, it follows that
\begin{equation}
  \label{eq:monotonicity-1}
\mathcal{S}_R \ge \mathcal{S}_{R'} \ge \mathcal{S}_\infty \qquad \text{for $R'>R>1$} \, .
\end{equation}


We claim that
\begin{equation}\label{min-infinite-strip-attained}
\mbox{the infimum in \eqref{sobolevconstants1-star} is attained at a function $v \in H^1_0(\Omega_\infty)$ with $\|v\|_{L^4(\Omega_\infty)}=1$.}
\end{equation}
In order to prove \eqref{min-infinite-strip-attained} we consider a minimizing sequence for \eqref{sobolevconstants1-star}, namely
a sequence $(u_n)_{n \in \mathbb{N}} \subset H^1_0(\Omega_\infty)$ such that $\|u_n\|_{L^4(\Omega_{\infty})}=1$, for every $n \in \N$,
and $\|\nabla u_n\|^2_{L^2(\Omega_\infty)} \to \mathcal{S}_{\infty}$. It then follows from the Poincaré inequality that $(u_n)_{n\in \mathbb{N}}$
is bounded in $H^1_0(\Omega_\infty)$. By Lions' Lemma \cite{lions1984concentration}, after passing to a subsequence, there exist points
$\tau_n=(x_n,y_n) \in \R^2$, $n \in \N$, with
$$
v_n\doteq \tau_n * u_n \rightharpoonup v \qquad \text{weakly in $H^1(\R^2)$ for some function $v \in H^1(\R^2) \setminus \{0\}$},
$$
in which we use the notation $\tau * u\doteq u(\,\cdot\, - \tau)$ to denote the translation of a function $u: \R^2 \longrightarrow \R$ with respect
to $\tau \in \R^2$. Since $\Omega_\infty$ is bounded in the $y$-direction and $v \not = 0$, also the sequence $(y_n)_{n \in \N}$ is bounded
and we may therefore assume that $y_n= 0$ for every $n \in \N$. Consequently, $v \in H^1_0(\Omega_\infty)$ and $v_n \rightharpoonup v$ in $H^1_0(\Omega_\infty)$. By the Brezis-Lieb Lemma \cite{brezis1983relation}, we then have
$$
1 = \lim_{n \to \infty}\|u_n\|^{4}_{L^4(\Omega_\infty)} = \lim_{n \to \infty}\|v_n\|^{4}_{L^4(\Omega_\infty)}
= \|v\|^{4}_{L^4(\Omega_\infty)} + c^4 \qquad \text{with}\quad c\doteq \lim_{n \to \infty}\|v_n-v\|_{L^4(\Omega_\infty)} \, .
$$
Consequently, $\|v\|^{2}_{L^4(\Omega_\infty)} + c^2 \ge 1$ and therefore
\begin{equation}
\mathcal{S}_{\infty} \le \mathcal{S}_{\infty} \left(\|v\|^{2}_{L^4(\Omega_\infty)} + c^2  \right) \le \|\nabla v\|^{2}_{L^2(\Omega_\infty)} + \lim_{n \to \infty}\|\nabla (v_n-v)\|^{2}_{L^2(\Omega_\infty)} = \lim_{n \to \infty}\|\nabla v_n\|^{2}_{L^2(\Omega_\infty)} = \mathcal{S}_{\infty}. \label{convexity-ineq}
\end{equation}
From this we deduce that
$$
\|v\|^{2}_{L^4(\Omega_\infty)} + c^2 = 1= \|v\|^{4}_{L^4(\Omega_\infty)} + c^4 \, .
$$
Since $v\neq0$, we infer that $c=0$ and then $\|v\|_{L^4(\Omega_\infty)}=1$. Therefore,
$\|\nabla v\|^{2}_{L^2(\Omega_\infty)} \ge \mathcal{S}_{\infty}$ by definition of $\mathcal{S}_{\infty}$ in
\eqref{sobolevconstants1-star}.
Then it follows from (\ref{convexity-ineq}) that $\|\nabla v\|^{2}_{L^2(\Omega_\infty)} = \mathcal{S}_\infty$ and $\lim \limits_{n \to \infty}\|\nabla (v_n-v)\|_{L^2(\Omega_\infty)} = 0$, and consequently $\lim \limits_{n \to \infty}\|v_n-v\|_{L^2(\Omega_\infty)} = 0$ by the Poincar\'e inequality. Hence $v_n \to v$ in $H^1_0(\Omega_\infty)$, and $v$ is a minimizer of the of the Sobolev quotient (\ref{sobolevconstants1-star}). This proves the claimed statement \eqref{min-infinite-strip-attained}.\par
We next claim that
\begin{equation}\label{eq:proof-asymptotic-1}
\lim_{R \to \infty} \mathcal{S}_{R}= \mathcal{S}_{\infty}\, .
\end{equation}
To show this, let $\rho \in \mathcal{C}^\infty_{0}(\R^2)$ be a nonnegative function with $\rho \equiv 1$ on $B_{\frac{1}{2}}$ and
$\rho \equiv 0$ on $\R^2 \setminus B_{1}$. Moreover, let
$$
\rho_n(x) \doteq \rho \left(\frac{x}{n} \right) \qquad \forall x \in \R^{2}\, , \ n \geq 1 ,
$$
so that $\rho_n \in \mathcal{C}^\infty_{0}(B_{n})$. By \eqref{min-infinite-strip-attained}, there exists a function $v \in H^1_0(\Omega_\infty)$ with $\|v\|_{L^4(\Omega_\infty)}=1$ and $\|\nabla v\|_{L^2(\Omega_\infty)}^2 = \mathcal{S}_{\infty}$.
It is then standard to see that the sequence $v_n\doteq \rho_n v \in H^1_0(\Omega_\infty)$, $n \in \N$, satisfies $v_n \to v$ in $H^1_0(\Omega_\infty)$ and, hence,
  \begin{equation}
    \label{eq:v-n-cnv}
  \|v_n\|_{L^4(\Omega_\infty)} \to 1 \quad \text{and}\quad \|\nabla v_n\|_{L^2(\Omega_\infty)}^2 \to \mathcal{S}_{\infty} \qquad \text{as} \ \ n \to \infty \, .
  \end{equation}
Since $v_n = 0$ on $\Omega_\infty \setminus B_n$ we have
  $$
  u_n \doteq x_n * v_n \in H^1_0(\Omega_{R_n})\quad \text{for $n \in \N$ with $x_n \doteq (n+1,0)$ and $R_n\doteq 2n+1$.}
  $$
  It thus follows from (\ref{eq:monotonicity-1}) and (\ref{eq:v-n-cnv}) that
  $$
  \mathcal{S}_{\infty} \le \lim_{R \to \infty} \mathcal{S}_{R} = \lim_{n \to \infty} \mathcal{S}_{R_n} \le \lim_{n \to \infty} \frac{\|\nabla u_n\|_{L^2(\Omega_{R_n})}^2}{\|u_n\|_{L^4(\Omega_{R_n})}^2}  = \lim_{n \to \infty} \frac{\|\nabla v_n\|_{L^2(\Omega_\infty)}^2}{\|v_n\|_{L^4(\Omega_\infty)}^2} = \mathcal{S}_{\infty}
  $$
  which yields the equality in (\ref{eq:proof-asymptotic-1}).\par
By \eqref{eq:proof-asymptotic-1}, the proof of the theorem will be finished after showing that
\begin{equation}
  \label{eq:proof-asymptotic-2}
\lim_{R \to \infty} \mathcal{S}_{R,e}> \mathcal{S}_{\infty}.
\end{equation}
To prove (\ref{eq:proof-asymptotic-2}), we argue by contradiction, assuming that there exists an increasing divergent sequence $(R_n)_{n \in \N} \subset (1,\infty)$ and, for every $n \in \N$, a function $u_n \in H^1_{0,e}(\Omega_{R_n})$ such that $\|u_n\|_{L^4(\Omega_{R_n})}=1$ and
$$
\|\nabla u_n\|^2_{L^2(\Omega_{R_n})} = \|\nabla u_n\|^2_{L^2(\Omega_\infty)} \to \mathcal{S}_{\infty} \qquad \text{as $n \to \infty$.}
$$
Without loss of generality, we may assume that $u_n \ge 0$ for every $n$, otherwise we replace $u_n$ by $|u_n|$. From the proof of \eqref{min-infinite-strip-attained}, we know that there exists a sequence $(t_n)_{n \in \mathbb{N}} \subset \R$ and a minimizer $v \in H^1_0(\Omega_\infty)$ of the Sobolev quotient (\ref{sobolevconstants1-star}) with $\|v\|_{L^4(\Omega_\infty)}=1$ and the property that
$$
(t_n,0) * u_n \to v \qquad \text{in $H^1_0(\Omega_\infty)$.}
$$
Since $u_n$ is symmetric with respect to the $x$-variable, we have
$$
[(-t_n,0) * u_n](x,y)= u_n(x+t_n,y)= u_n(-x-t_n,y) =[(t_n,0)* u_n](-x,y) \qquad \forall (x,y) \in \Omega_{R_n} \, , \ n \in \N \, ,
$$
which implies that
$$
(-t_n,0) * u_n \to \tilde v \qquad \text{in $H^1_0(\Omega_\infty)$,}
$$
where $\tilde v \in H^1_0(\Omega_\infty)$ is defined by $\tilde v(x,y) \doteq v(-x,y)$, for every $(x,y) \in \Omega_{R}$. Then, we
note the following elementary fact:
if, for a sequence $(w_n)_{n \in \mathbb{N}} \subset H^1(\R^2)$ and some sequences $(\alpha_n)_{n \in \mathbb{N}}, (\beta_n)_{n \in \mathbb{N}} \subset \R^2$ we have that both $(\alpha_n * w_n)_{n \in \mathbb{N}}$ and $(\beta_n * w_n)_{n \in \mathbb{N}}$ are convergent sequences in $H^1(\R^2)$, then the sequence $(\alpha_n-\beta_n)_{n \in \mathbb{N}}$ is bounded in $\R^2$.
From this fact, it follows that $(t_n)_{n \in \N}$ is bounded. Passing to a subsequence we may then assume, without loss of generality, that $t_n = 0$ for every $n \in \N$. Hence,
$$
u_n \to v \qquad \text{in $H^1_0(\Omega_\infty)$ as $n \to \infty$.}
$$
Since $v$ is a minimizer of the Sobolev quotient (\ref{sobolevconstants1-star}) with $\|v\|_{L^4(\Omega_\infty)}=1$, it follows that $v$ is a nontrivial weak solution of the semilinear elliptic equation (see, in particular, \cite{amick1984semilinear})
\begin{equation}\label{deq-limit}
\left\{
\begin{aligned}
& -\Delta v = \mathcal{S}_{\infty} v^3\ &&\qquad \text{in} \quad \Omega_{\infty} \, , \\[5pt]
& v = 0 &&\qquad \text{on} \quad \partial\Omega_{\infty} \, .
\end{aligned}
\right.
\end{equation}
By standard elliptic regularity we have that $v \in \mathcal{C}^\infty(\Omega_\infty)$. Moreover, $v \ge 0$ in $\Omega_\infty$ and $v = 0$ in $B_1$ since $u_n \in H^1_0(\Omega_{R_n})$ for every $n \in \N$. Since $v$ is superharmonic in $\Omega_\infty$ by (\ref{deq-limit}), this contradicts the strong maximum principle. The proof of (\ref{eq:proof-asymptotic-2}) is thus finished, which completes the proof of the theorem.\end{proof}

\begin{remark} \label{conn1}
Theorem \ref{tobias} establishes a first connection between symmetry breaking in the stationary Navier-Stokes equations \eqref{nsstokes0} and Sobolev minimizers \eqref{sobolevconstants1}. In fact, the spirit of its proof is that the mass of the Sobolev minimizer in \eqref{sobolevconstants1} tends to concentrate on one side of the obstacle as $R \to \infty$. This is similar to the phenomenon depicted in Figure \ref{vortexshedding}, where the energy (vortex shedding) mainly concentrates behind the obstacle as the Reynolds number increases.
\end{remark}

As a straightforward consequence of Remark \ref{remrec} and \eqref{eq:proof-asymptotic-2} we obtain the estimate
\begin{equation} \label{bounds1}
\dfrac{\pi}{2} \sqrt{ \dfrac{3 \pi}{2} } \dfrac{1}{h} \leq \mathcal{S}_{\infty} \leq \dfrac{1}{2} \sqrt{\dfrac{\pi}{2 \kappa}} (2 h^2 \log(h)^2 - 2 h^2 \log(h) + h^2 - 1) \qquad \forall h >1 \, .
\end{equation}
A finer upper bound for $\mathcal{S}_{\infty}=\mathcal{S}_{\infty}(h)$ can be found by seeking the minimum in \eqref{sobolevconstants1-star} among separated-variable functions, that is, having the form $V(x) W(y)$ for $(x,y) \in \Omega_{\infty}$.
\begin{proposition} \label{betterbound}
For every $h > 1$ we have
\begin{equation} \label{bounds2}
\dfrac{\pi}{2} \sqrt{ \dfrac{3 \pi}{2} } \dfrac{1}{h} \leq \mathcal{S}_{\infty}(h) \leq \dfrac{5.151}{h} \, .
\end{equation}	
\end{proposition}
\noindent
\begin{proof}
Given $h>1$ and $W_{h} \in H^{1}_{0}(-h,h) \setminus \{0\}$, \eqref{min-infinite-strip-attained} implies that
\begin{equation} \label{sepvar1}
\mathcal{S}_{\infty}(h) \leq \min_{V \in H^{1}(\R) \setminus \{0\}} \ \frac{\|V' \|^{2}_{L^2(\R)} \|W_{h} \|^{2}_{L^2(-h,h)} + \|V \|^{2}_{L^2(\R)} \|W'_{h} \|^{2}_{L^2(-h,h)}}{\|V \|^{2}_{L^4(\R)} \|W_{h} \|^{2}_{L^4(-h,h)}} \, .
\end{equation}
Furthermore, a rescaling argument shows that
\begin{equation} \label{sepvar2}
\mathcal{S}_{\infty}(h) = \dfrac{h_{0}}{h} \mathcal{S}_{\infty}(h_{0}) \qquad \forall h, h_{0} > 1 \, .
\end{equation}
Let $ \text{cn} : \mathbb{R} \longrightarrow [-1,1]$ be the Jacobian elliptic cosine function with modulus $k = 1/ \sqrt{2}$, which satisfies the one-dimensional version of \eqref{deq-limit}, namely
$$
\text{cn}''(t) + \text{cn}(t)^3 = 0 \qquad \forall t \in \mathbb{R} \, ,
$$
and whose first zero is given by
$$
\alpha \doteq \sqrt{2} \int_{0}^{\pi/2} \dfrac{1}{\sqrt{2 - \sin(t)^2}} \, dt \approx 1.85407,
$$
see \cite{armitage2006elliptic} for more details. Then, the function
$$
W_{h}(y) \doteq \dfrac{1}{\mu_{h}} \text{cn} \left( \dfrac{\alpha}{h} y \right) \qquad \forall y \in [-h,h]
$$
vanishes at $y=\pm h$, where $\mu_{h}>0$ is a normalization constant such that $\|W_{h} \|_{L^4(-h,h)}=1$. Let $h_{0} > 1$ be such that  $\|W'_{h_0} \|_{L^2(-h_{0},h_{0})}=1$; numerically we find $h_{0} \approx 1.98978$. Then, the Euler-Lagrange equation associated to the minimization problem in \eqref{sepvar1} (when $h=h_0$) reads
\begin{equation} \label{sepvar3}
- \|W_{h_{0}} \|^{2}_{L^2(-h_{0},h_{0})} V''(x) + V(x) = \lambda V(x)^3 \qquad \forall x \in \R \, ,
\end{equation}
with $\lambda \in \R$ being a Lagrange multiplier. By direct substitution we deduce that the function
$$
V_{h_{0}}(x) \doteq \left[ \cosh \left( \dfrac{x}{\|W_{h_{0}} \|_{L^2(-h_{0},h_{0})}} \right) \right]^{-1}  \qquad \forall x \in \R \,
$$
is a solution of \eqref{sepvar3} with $\lambda =2$. The function $V_{h_{0}}W_{h_{0}}$ minimizes the ratio in \eqref{sepvar1} when $W_{h}=W_{h_0}$, and its graph is depicted in Figure \ref{jacobi} below. 

\begin{figure}[H]
	\begin{center}
		\includegraphics[height=65mm,width=90mm]{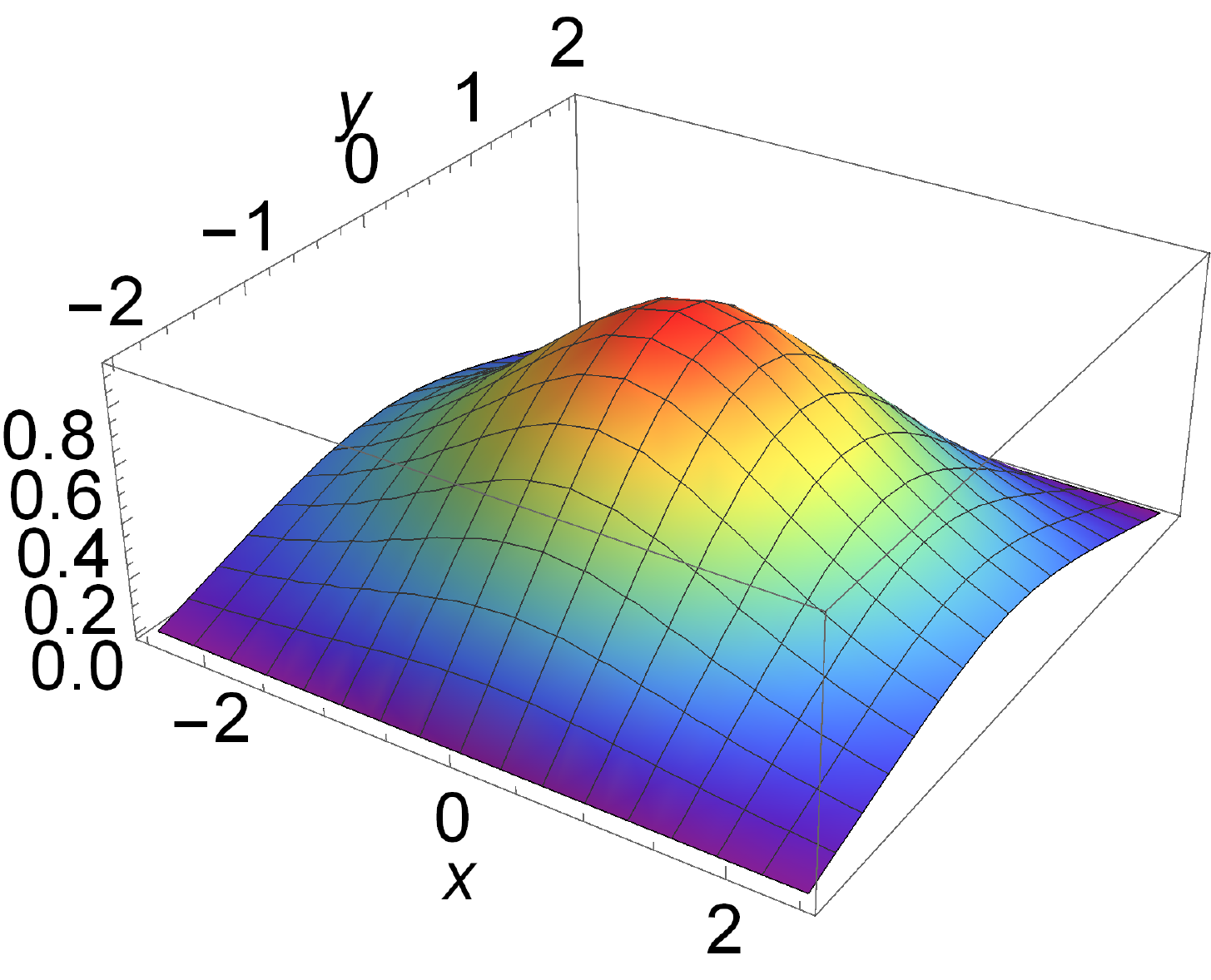}
	\end{center}
	\vspace*{-5mm}
	\caption{Graph of the function $V_{h_{0}}(x)W_{h_{0}}(y)$ for $(x,y) \in [-2.5,2.5] \times [-h_0, h_0]$.}\label{jacobi}
\end{figure}

The approximate value $5.151$ is obtained numerically after replacing the functions $V_{h_{0}}$ and $W_{h_{0}}$ into the ratio \eqref{sepvar1} (with $h=h_{0}$) and using \eqref{sepvar2}.
\end{proof}

The ratio between the upper and lower bounds for $\mathcal{S}_{\infty}$ found in Proposition \ref{betterbound} is approximately $1.51061$, showing the accuracy of these bounds. A comparison between the lower and upper bounds in \eqref{bounds1}-\eqref{bounds2} is shown in Figure \ref{jacobi2} below.
\begin{figure}[H]
	\begin{center}
		\includegraphics[height=70mm,width=160mm]{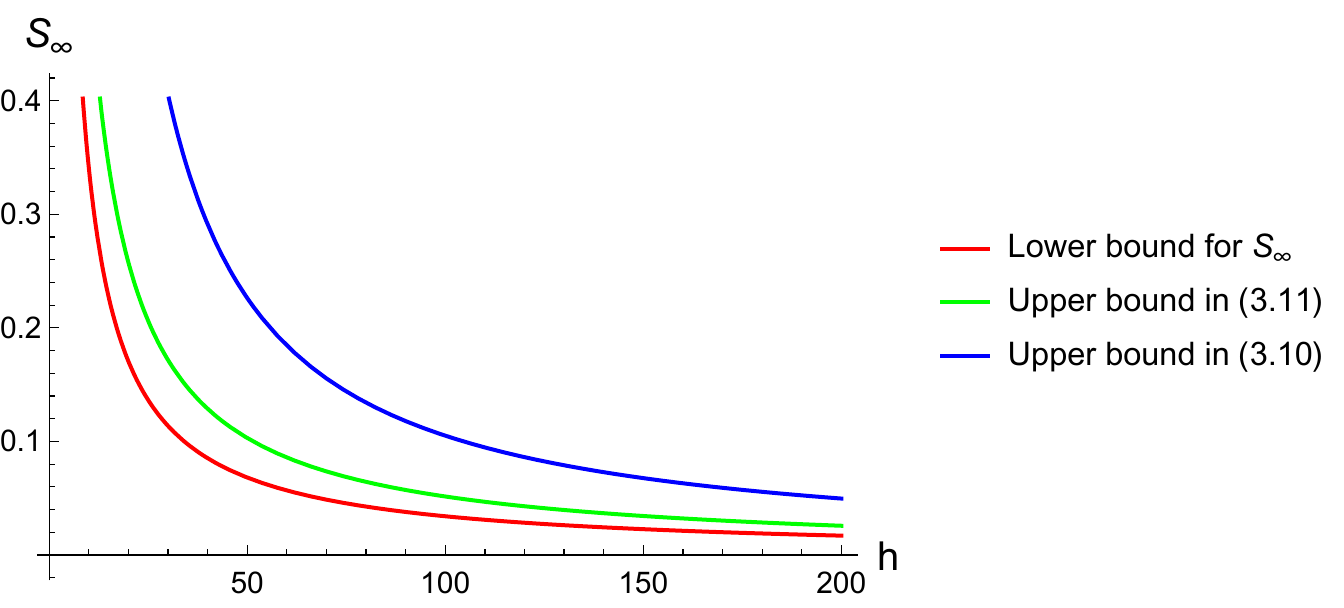}
	\end{center}
	\vspace*{-5mm}
	\caption{Comparison between the lower and upper bounds for $\mathcal{S}_{\infty} = \mathcal{S}_{\infty}(h)$, for $h >1$.}
	\label{jacobi2}
\end{figure}

\newpage
\section{Symmetry breaking in the Navier-Stokes equations}

\subsection{Non-uniqueness as a necessary condition for symmetry breaking}

In this subsection we start by quickly recalling some well-known functional spaces and inequalities, adapting them to our context.
Let us introduce the two functional spaces of vector fields
$$
\mathcal{V}_{c}(\Omega_{R}) = \{ v \in H^1_c(\Omega_{R})^2 \ | \ \nabla \cdot v = 0 \ \ \text{in} \ \Omega_{R} \} \qquad \text{and} \qquad \mathcal{V}(\Omega_{R}) = \{ v \in H^{1}_{0}(\Omega_{R}) \ | \ \nabla \cdot v = 0 \ \ \text{in} \ \Omega_{R} \},
$$
which are Hilbert spaces if endowed with the scalar product $(u,v)\mapsto(\nabla u,\nabla v)_{L^2(\Omega_{R})}$. We also introduce the continuous trilinear form
$$
\beta(u,v,w) = \int_{\Omega_{R}} (u\cdot\nabla)v \cdot w \qquad \forall u, v, w \in H^{1}(\Omega_{R}),
$$
which satisfies, by H\"older's inequality,
\begin{equation} \label{propiedadestri0}
| \beta(u,v,w) | \leq \| u \|_{L^{4}(\Omega_{R})} \| \nabla v \|_{L^{2}(\Omega_{R})} \| w \|_{L^{4}(\Omega_{R})} \quad\qquad \forall u,v,w \in H^1(\Omega_{R}) \, .
\end{equation}
Moreover,
\begin{equation} \label{propiedadestri}
\begin{array}{ll}
\beta(u,v,w) = - \beta(u,w,v) \ \ \ \ &\text{for any} \ \ \ u \in \mathcal{V}_{c}(\Omega_{R}), \ v \in H^{1}(\Omega_{R}), \ w \in H_{0}^{1}(\Omega_{R}),\\[3pt]
\beta(u,v,v) = 0 \ \ \ \ &\text{for any} \ \ \ u \in \mathcal{V}_{c}(\Omega_{R}), \ v \in H_{0}^{1}(\Omega_{R}).
\end{array}
\end{equation}

We are now in position to define weak solutions for problem \eqref{nsstokes0}.

\begin{definition}\label{weaksolution}
Given $U>0$, we say that a vector field $u \in \mathcal{V}_{c}(\Omega_{R})$ is a weak solution of \eqref{nsstokes0} if $u$ verifies \eqref{nsstokes0}$_2$ in the trace sense and
	\begin{equation} \label{nstokesdebil}
	\eta (\nabla u, \nabla \varphi)_{L^{2}(\Omega_{R})} + \beta(u,u,\varphi) = 0 \qquad \forall \varphi \in \mathcal{V}(\Omega_{R}).
	\end{equation}	
\end{definition}

The existence of weak solutions of \eqref{nsstokes0} relies on the construction of a solenoidal extension of the boundary velocity, see \cite{ladyzhenskaya1978some}. Namely, one needs to find a vector field $\Psi\in H^{1}(\Omega_{R})$ such that
\begin{equation}\label{div}
\nabla\cdot\Psi=0\quad\mbox{in }\Omega_{R}\, ,\qquad \Psi=(U,0)\quad\text{on } \partial Q_{R} \, ,\qquad\Psi=(0,0)\quad\text{on }\partial B_{1} \, .
\end{equation}

Then we prove a result whose qualitative statement is essentially known, see for example \cite[Section IX.4]{galdi2011introduction}. We give here a slightly different statement and a full proof because explicit constants are needed to show how uniqueness depends on the magnitude of the solenoidal extension through the Sobolev embedding constant \eqref{sobolevconstants1}.

\begin{theorem}\label{UVconstant0}
For any $U>0$ there exists a weak solution $(u,p)\in\mathcal{V}_{c}(\Omega_{R}) \times L^{2}(\Omega_{R})$ of \eqref{nsstokes0}. If, moreover, there exists a vector field $\Psi \in H^{1}(\Omega_{R})$ verifying \eqref{div} and the inequality
\begin{equation} \label{umbral1}
2\| \nabla \Psi \|_{L^{2}(\Omega_{R})} +\sqrt{\mathcal{S}_{R}} \, \|\Psi\|_{L^4(\Omega_{R})} < \eta \mathcal{S}_{R} \, ,
\end{equation}
then the weak solution of \eqref{nsstokes0} is unique, satisfies the symmetry property
\begin{equation} \label{simpro}
u_1(x,-y)=u_1(x,y),\quad u_2(x,-y)=-u_2(x,y),\quad p(x,-y)=p(x,y)\qquad \forall (x,y) \in \Omega_{R} \, .
\end{equation}
and the estimate
\begin{equation} \label{estima}
\| \nabla u \|_{L^{2}(\Omega_{R})}  < 3 \| \nabla \Psi \|_{L^{2}(\Omega_{R})} < \dfrac{3}{2} \eta \mathcal{S}_{R}  \, .
\end{equation}
\end{theorem}
\noindent
\begin{proof}
Existence of a weak solution $(u,p)\in\mathcal{V}_{c}(\Omega_{R}) \times L^{2}(\Omega_{R})$ of \eqref{nsstokes0} follows from \cite[Theorem 3.1]{gazspe}, recalling that the compatibility condition \eqref{zeroflux2} is fulfilled.\par
In order to prove the uniqueness statement, suppose there exists a vector field $\Psi \in H^{1}(\Omega_{R})$ verifying \eqref{div} and \eqref{umbral1} and, by contradiction, that there exist two weak solutions $u,v\in \mathcal{V}_{c}(\Omega_{R})$ of \eqref{nsstokes0}. Set $\xi \doteq v-\Psi$, so that $\xi\in\mathcal{V}(\Omega_{R})$ can be used as a test function in the weak formulation \eqref{nstokesdebil} satisfied by $v$, thus yielding
$$
\eta\|\nabla\xi\|^{2}_{L^{2}(\Omega_{R})} + \eta (\nabla \Psi, \nabla \xi) + \beta(\xi + \Psi,\xi + \Psi,\xi) = 0
$$
and, in turn,
\begin{equation}\label{unicidad5}
\eta\|\nabla\xi\|^{2}_{L^{2}(\Omega_{R})} \leq \eta \| \nabla \Psi \|_{L^{2}(\Omega_{R})} \|\nabla \xi \|_{L^{2}(\Omega_{R})}
- \beta(\xi + \Psi,\xi + \Psi,\xi) .
\end{equation}
In view of \eqref{propiedadestri0}-\eqref{propiedadestri} we have $\beta(\xi + \Psi,\xi + \Psi,\xi) = \beta(\xi + \Psi,\Psi,\xi)$ and the
estimate
\begin{equation} \label{unicidad6}
\begin{aligned}
|\beta(\xi+\Psi,\Psi,\xi)| & \leq \|\xi\|_{L^4(\Omega_{R})}\|\nabla \Psi\|_{L^2(\Omega_{R})} \left( \|\xi\|_{L^4(\Omega_{R})}+\|\Psi\|_{L^4(\Omega_{R})} \right) \\[6pt]
& \leq \dfrac{\| \nabla \xi \|_{L^{2}(\Omega_{R})}}{\sqrt{\mathcal{S}_{R}}} \| \nabla \Psi \|_{L^{2}(\Omega_{R})}
\left( \dfrac{\| \nabla \xi \|_{L^{2}(\Omega_{R})}}{\sqrt{\mathcal{S}_{R}}} +\|\Psi\|_{L^4(\Omega_{R})}\right) \, ,
\end{aligned}
\end{equation}
where we used \eq{sobolevconstants11}. After inserting \eqref{unicidad6} into \eqref{unicidad5} and applying \eqref{umbral1} so that
$\| \nabla \Psi \|_{L^{2}(\Omega_{R})} < \mathcal{S}_{R} \eta $, we infer the bound
\begin{equation}\label{unicidad7}
\| \nabla \xi \|_{L^{2}(\Omega_{R})} \leq \dfrac{ \dfrac{\|\nabla \Psi\|_{L^2(\Omega_{R})}\|\Psi\|_{L^4(\Omega_{R})}}{\sqrt{\mathcal{S}_{R}}}+\eta\|\nabla \Psi \|_{L^{2}(\Omega_{R})}}{\eta - \dfrac{\| \nabla \Psi \|_{L^{2}(\Omega_{R})}}{\mathcal{S}_{R}}}.
\end{equation}
Let $w \doteq u - v \in \mathcal{V}(\Omega_{R})$ and
subtract the equations \eqref{nstokesdebil} corresponding to $u$ and $v$ in order to obtain
$$
\eta (\nabla w,\nabla\varphi)_{L^2(\Omega_{R})} + \beta(u,w,\varphi) + \beta(w,v,\varphi) = 0 \ \ \ \ \ \ \forall \varphi \in \mathcal{V}(\Omega_{R}).
$$
By taking $\varphi=w$ and using \eqref{propiedadestri0} and \eqref{unicidad7}, we derive
$$
\begin{aligned}
\eta\|\nabla w\|^2_{L^2(\Omega_{R})} & = - \beta(w,v,w)=\beta(w,w,v) \leq \| w\|_{L^{4}(\Omega_{R})} \| \nabla w \|_{L^{2}(\Omega_{R})}\| v\|_{L^{4}(\Omega_{R})} \\[6pt]
& \leq \dfrac{\| \nabla w \|^{2}_{L^{2}(\Omega_{R})}}{\sqrt{\mathcal{S}_{R}}}\, \|v \|_{L^4(\Omega_{R})} \leq \dfrac{\| \nabla w \|^{2}_{L^{2}(\Omega_{R})}}{\sqrt{\mathcal{S}_{R}}} \left(\|\xi\|_{L^4(\Omega_{R})}+\|\Psi\|_{L^4(\Omega_{R})} \right) \\[6pt]
& \leq \dfrac{\| \nabla w \|^{2}_{L^{2}(\Omega_{R})}}{\sqrt{\mathcal{S}_{R}}} \left(\dfrac{\| \nabla \xi \|_{L^{2}(\Omega_{R})}}{\sqrt{\mathcal{S}_{R}}} +\|\Psi\|_{L^4(\Omega_{R})} \right) \\[6pt]
& \leq \eta \| \nabla w \|^{2}_{L^{2}(\Omega_{R})} \, \dfrac{\| \nabla \Psi \|_{L^{2}(\Omega_{R})} +\sqrt{\mathcal{S}_{R}} \, \|\Psi\|_{L^4(\Omega_{R})}}{\eta\mathcal{S}_{R}-\|\nabla \Psi\|_{L^{2}(\Omega_{R})}} ,
\end{aligned}
$$
which shows that $w=0$ (and, therefore, unique weak solvability for \eqref{nsstokes0}) provided that \eqref{umbral1} holds. The estimate \eqref{estima} follows from \eqref{umbral1} and \eqref{unicidad7}, noticing that
$$
\begin{aligned}
\|\nabla u\|_{L^2(\Omega_{R})} = \|\nabla v\|_{L^2(\Omega_{R})} & \leq \|\nabla \xi\|_{L^2(\Omega_{R})} + \|\nabla \Psi\|_{L^2(\Omega_{R})} \\[6pt]
& \leq \dfrac{ \sqrt{\mathcal{S}_{R}} \, \|\nabla \Psi\|_{L^2(\Omega_{R})}\|\Psi\|_{L^4(\Omega_{R})} + 2 \eta \mathcal{S}_{R} \|\nabla \Psi \|_{L^{2}(\Omega_{R})} - \|\nabla \Psi\|^{2}_{L^2(\Omega_{R})}}{\eta \mathcal{S}_{R} -\| \nabla \Psi \|_{L^{2}(\Omega_{R})}} \\[6pt]
& < \dfrac{  \|\nabla \Psi\|_{L^2(\Omega_{R})} \left( \eta \mathcal{S}_{R} - 2\| \nabla \Psi \|_{L^{2}(\Omega_{R})} \right) + 2 \eta \mathcal{S}_{R}  \|\nabla \Psi \|_{L^{2}(\Omega_{R})} - \|\nabla \Psi\|^{2}_{L^2(\Omega_{R})}}{\eta \mathcal{S}_{R} -\| \nabla \Psi \|_{L^{2}(\Omega_{R})}} \\[6pt]
& = 3 \| \nabla \Psi \|_{L^{2}(\Omega_{R})} < \dfrac{3}{2} \eta \mathcal{S}_{R} \, .
\end{aligned}
$$

In case of uniqueness, the symmetry property \eqref{simpro} is established by noticing that the pair $(v_1,v_2,q) \in H^{1}(\Omega_{R}) \times L^{2}(\Omega_{R})$ defined as
$$
v_1(x,y)=u_1(x,-y),\quad v_2(x,y)=-u_2(x,-y),\quad q(x,y)=p(x,-y)\qquad\mbox{for a.e. }(x,y) \in \Omega_{R} \, ,
$$
also solves \eqref{nsstokes0}.\end{proof}

As a consequence of Theorem \ref{UVconstant0}, we infer that symmetry breaking for the flow, see Figure \ref{vortexshedding}, may occur only if
\eq{umbral1} is violated, that is, if {\em for any} vector field $\Psi \in H^{1}(\Omega_{R})$ verifying \eqref{div} one has
$$
2\| \nabla \Psi \|_{L^{2}(\Omega_{R})} +\sqrt{\mathcal{S}_{R}} \, \|\Psi\|_{L^4(\Omega_{R})} \geq \mathcal{S}_{R} \eta \, .
$$

In the next subsection we give a sufficient condition for the existence of a vector field $\Psi\in H^{1}(\Omega_{R})$ satisfying \eqref{div}-\eqref{umbral1}.

\subsection{A small solenoidal extension of the boundary velocity} \label{solextensions}

The purpose of this subsection is to build an explicit solenoidal extension of the boundary velocity, see \cite{ladyzhenskaya1978some}. We will build a vector field $\Psi\in H^{1}(\Omega_{R})$ satisfying \eqref{div} and whose ``size" can be explicitly computed or estimated in terms of $R$ and $h$. By exploiting the particular geometry of $\Omega_R$, a solenoidal extension is built through the use of a suitable cut-off function. Since we are interested in obtaining explicit bounds, the well-known Hopf cut-off function \cite[Lemma 1.2]{galdi1991existence} does not serve our purpose.

\begin{theorem} \label{rectext}
	Let $U > 0$ and define the constants
	$$
	\mathcal{B}_{1} = \dfrac{8}{35} \left[ \dfrac{36}{5} \dfrac{2h^2 + 3h +2}{(h-1)^2} + \dfrac{6}{5} \dfrac{(13h + 22)(3h^2 - h + 3)}{(h-1)^3} + \dfrac{19h^3 + 51h^2 + 75h +65}{3(h-1)^3} \right] \, ,
	$$
	\vspace{0.2cm}
	$$
	\begin{aligned}
	\mathcal{B}_{2} & = 4Rh + \dfrac{4 (82563626 + 139273674 h + 131633079  h^2 + 47395086 h^3)}{75150075 (h-1)} \\[3pt]
	& \hspace{4mm} + \dfrac{4 (6562533 + 20038773 h + 29176308 h^2 + 22648263 h^3 + 5977793 h^4)}{25050025  (h -1)^2} \\[3pt]
	& \hspace{4mm} + \dfrac{288 (1561958 + 3280874 h + 4160951 h^2 + 3491837 h^3 + 1768313 h^4 +
		336653 h^5)}{425850425  (h-1)^3} \, .
	\end{aligned}
	$$
	For every $R > 1$ there exists a vector field $\Psi_{R} \in H^{1}(\Omega_{R})$ satisfying \eqref{div} and
	\begin{equation} \label{vecpsi3}
	\|\nabla\Psi_{R}\|_{L^{2}(\Omega_{R})}=\sqrt{\mathcal{B}_{1}}\, U\ ,\qquad\|\Psi_{R}\|_{L^{4}(\Omega_{R})}=\sqrt[4]{\mathcal{B}_{2}}\, U\, .
	\end{equation}
\end{theorem}
\noindent
\begin{proof}
	Given any $\varepsilon > 0$ we define the function $\phi_{\varepsilon} : \mathbb{R} \longrightarrow \mathbb{R}$ as
	$$
	\phi_{\varepsilon} (t) = \begin {cases}
	0  & \text{ if } \ \ t \in (-\infty,-1-\varepsilon] \cup [1+\varepsilon,\infty)
	\\ \noalign{\medskip}
	\dfrac{1}{\varepsilon^{3}} \left[ 2 |t|^3 - 3(\varepsilon+2)t^{2} + 6(1+\varepsilon)|t| + \varepsilon^{3} - 3\varepsilon - 2 \right]  & \text { if } \ \ t \in (-1-\varepsilon,-1) \cup (1,1+\varepsilon)
	\\ \noalign{\medskip}
	1 & \text{ if } \ \ t \in [-1,1] \, ,
	\end{cases}
	$$
	see also \cite[Section 4]{sperone2021steady}, whose plot for $\varepsilon = 1/2$ is displayed in Figure \ref{sym3}. Then we have $\phi_{\varepsilon} \in \mathcal{C}^{1}(\mathbb{R})$, $\text{supp}(\phi_{\varepsilon}) =[-1-\varepsilon,1+\varepsilon]$, $\phi_{\varepsilon}'(1)=\phi_{\varepsilon}'(-1)=0$, so that  $\phi_{\varepsilon} \in H^{2}(\mathbb{R})$. In fact, it holds $\phi_{\varepsilon} \in W^{2,\infty}(\mathbb{R})$ with
	$$
	\| \phi_{\varepsilon} \|_{L^{\infty}(\mathbb{R})} =1 \, , \qquad
	\| \phi_{\varepsilon}' \|_{L^{\infty}(\mathbb{R})} = \dfrac{3}{2\varepsilon} \, , \qquad \| \phi_{\varepsilon}'' \|_{L^{\infty}(\mathbb{R})} = \dfrac{6}{\varepsilon^{2}} \, .
	$$	
	We then notice that
	\begin{equation} \label{propsigma2}
	\phi_{h-1}(\pm h)=0 \, ; \qquad \hspace{1mm} \phi_{h-1}(y)=1 \qquad \forall y \in [-1,1] \, ; \qquad \phi_{h-1}'(\pm 1)=\phi_{h-1}'(\pm h) = 0 \, .
	\end{equation}
	
		\begin{figure}[H]
		\begin{center}
			\includegraphics[height=40mm,width=80mm]{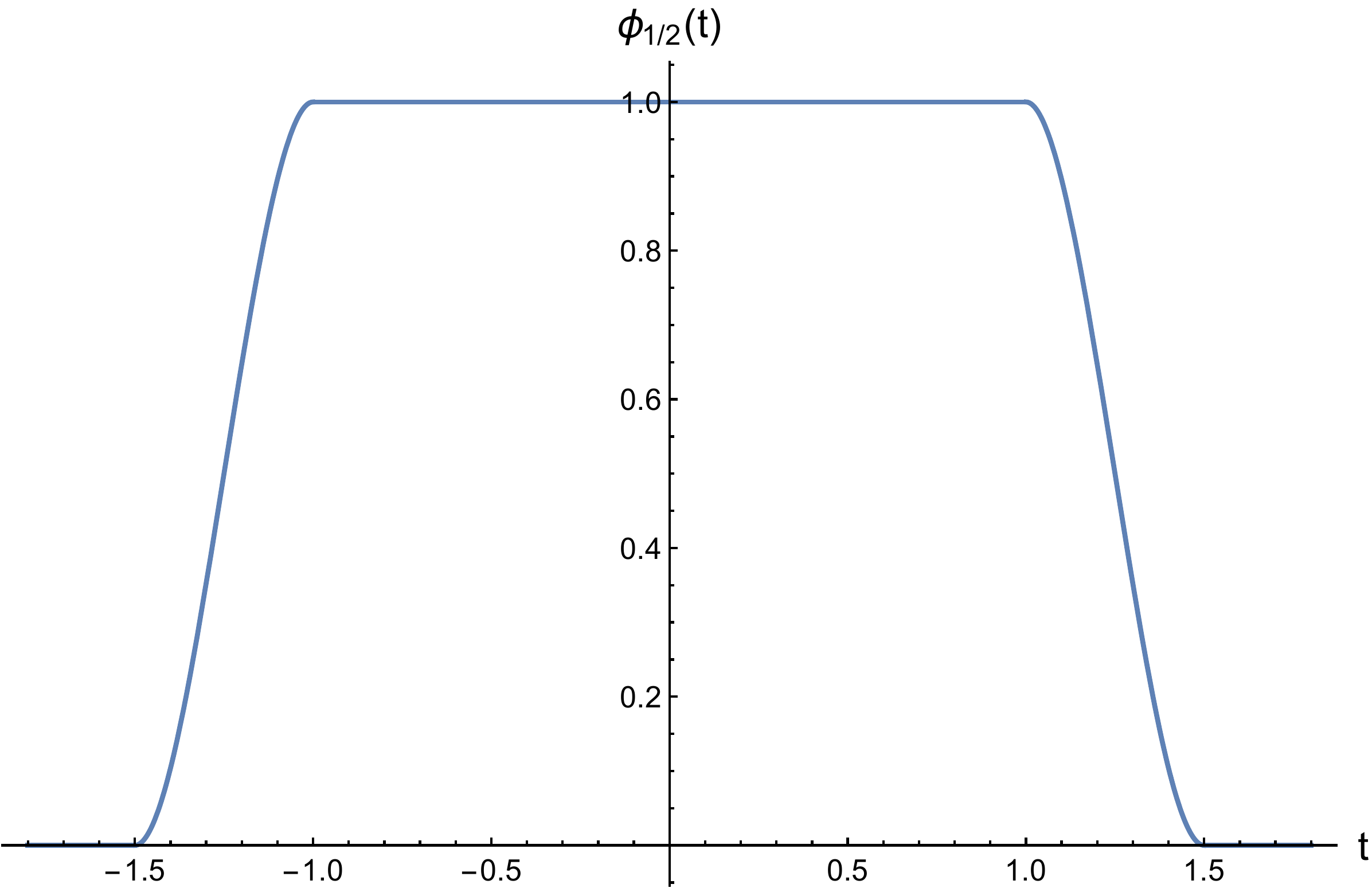} \quad \includegraphics[height=40mm,width=80mm]{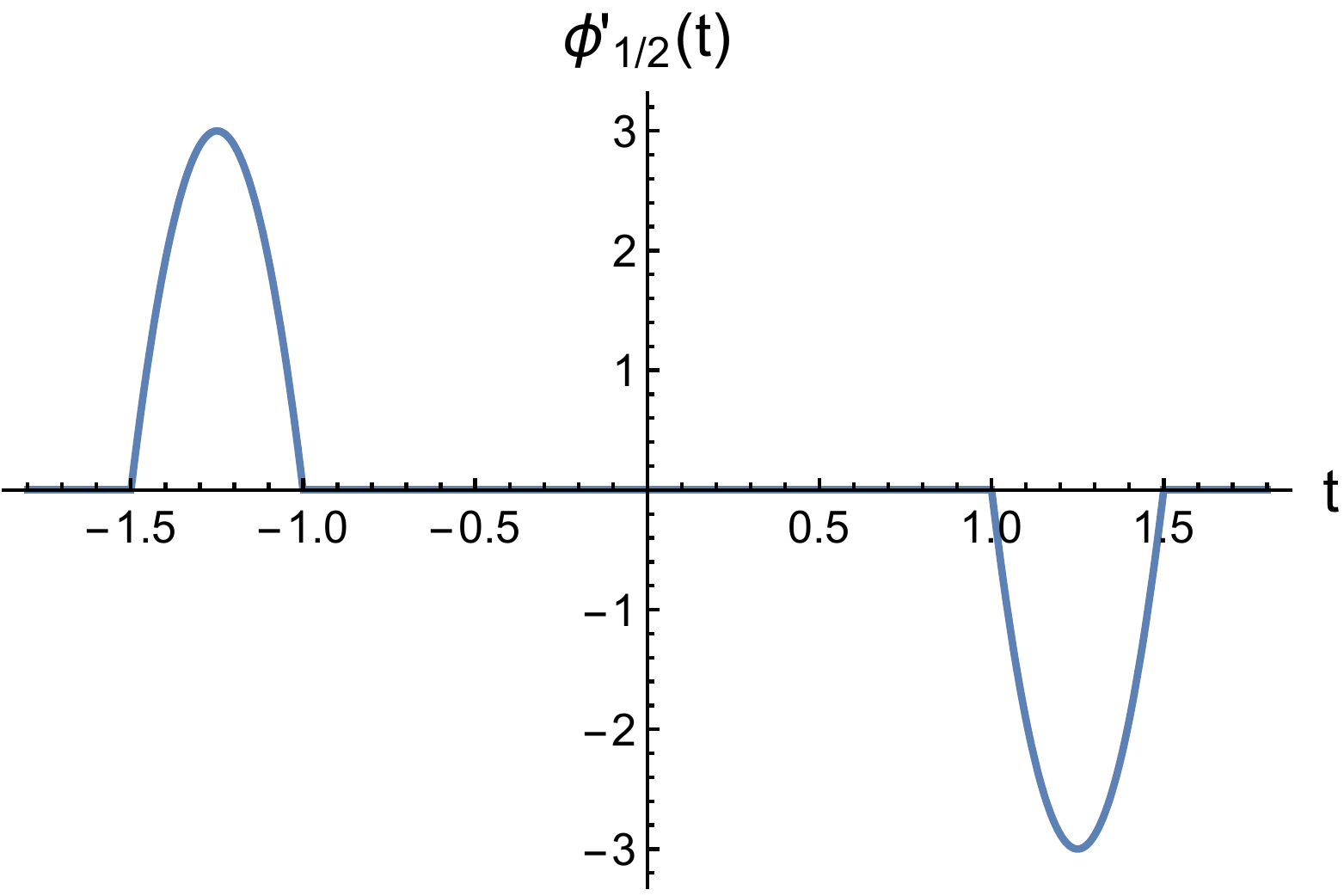}
		\end{center}
		\vspace*{-3mm}
		\caption{Graph of the function $\phi_{\varepsilon}$ (left) and of its derivative (right), for $\varepsilon = 1/2$.}\label{sym3}
	\end{figure}
	
	We then define the function $\omega \in H^{2}(Q_{R})$ by
	$$
	\omega(x,y) = 1 - \phi_{h-1}(x) \phi_{h-1}(y) \qquad \forall (x,y) \in Q_R \ ,
	$$
	whose plot is given in Figure \ref{cutoff3}. Notice that $\omega(x,y)=1$ when $|x| \geq h$, for every $y \in [-h,h]$.
	\begin{figure}[H]
		\begin{center}
			\includegraphics[height=65mm,width=90mm]{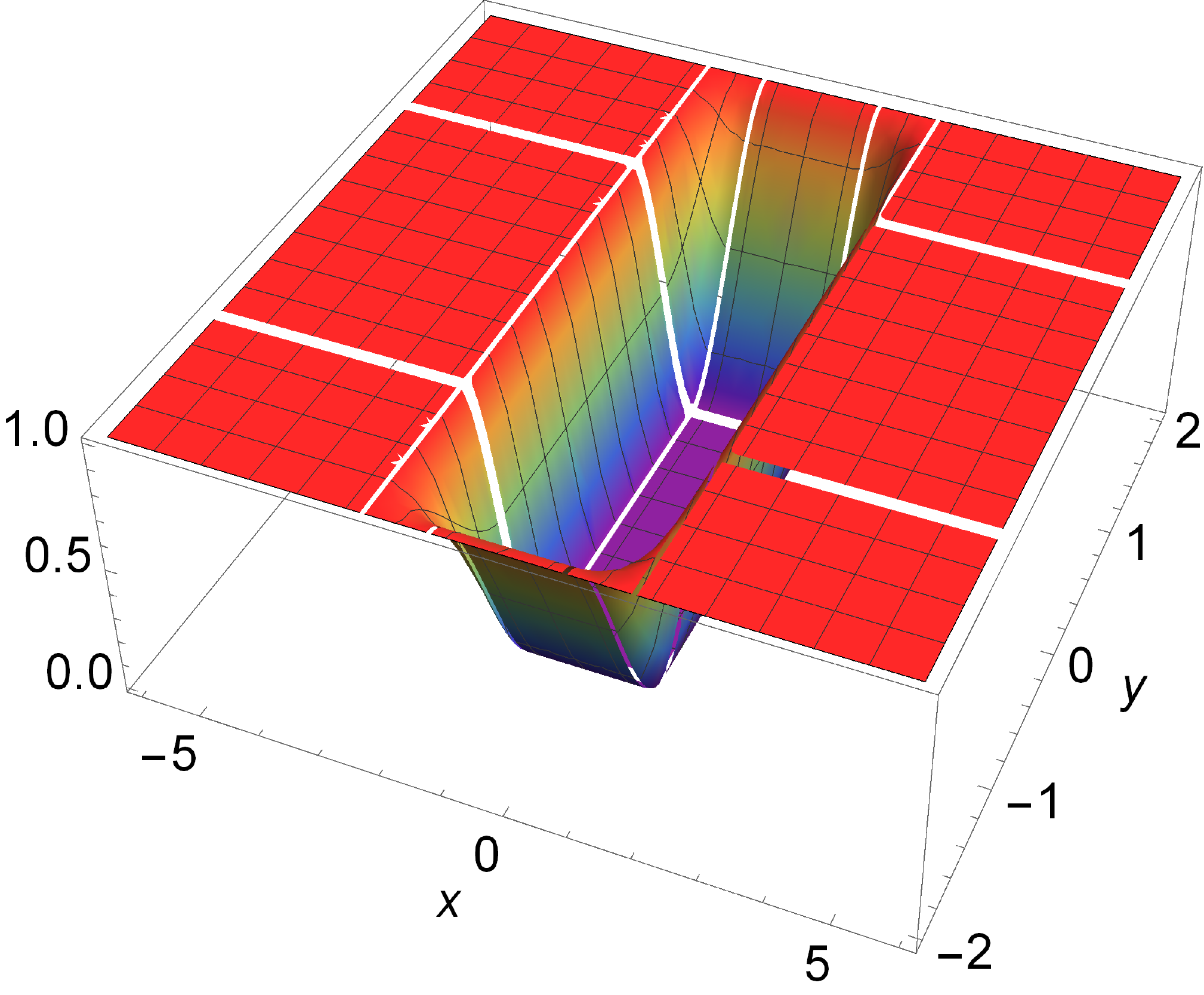}
		\end{center}
		\vspace*{-5mm}
		\caption{Graph of the function $\omega$ for $R =6$ and $h=2$.}\label{cutoff3}
	\end{figure}
	
	We are now in position to define the vector field
	$$
	\Psi_{R}(x,y) = U
	\left(
	\begin{matrix}
	\omega(x,y) + y \dfrac{\partial \omega}{\partial y}(x,y) \\[12pt]
	-y \dfrac{\partial \omega}{\partial x}(x,y)
	\end{matrix}
	\right)
	\qquad \forall (x,y) \in Q_R \, ,
	$$
	which is an element of $H^{1}(\Omega_{R})$ satisfying \eqref{div} in view of \eqref{propsigma2}. We have $\Psi(x,y)=(U,0)$ when $|x| \geq h$, for every $y \in [-h,h]$. The estimates \eqref{vecpsi3} follow from an explicit (but tedious) computation.
\end{proof}

The following result can be derived as a straightforward consequence of Theorems \ref{boundsobolevrec}, \ref{UVconstant0}  and \ref{rectext}.

\begin{corollary} \label{UVconstantcor}
	For any $U>0$ there exists a weak solution of \eqref{nsstokes0}. If, moreover,
	$$
	\dfrac{U}{\eta} <  \dfrac{\dfrac{\pi}{2} \sqrt{\dfrac{3 \pi}{2}} \dfrac{\sqrt{R^2 + h^2}}{Rh} } {2 \sqrt{\mathcal{B}_{1} } + \sqrt{\dfrac{\pi}{2Rh}} \, \sqrt[4]{ \dfrac{3 \pi \mathcal{B}_{2}}{2}(R^2 + h^2)} } \, ,
	$$
	then the weak solution of \eqref{nsstokes0} is unique and satisfies the symmetry property \eqref{simpro}.
\end{corollary}

Finally, we show that the family of solenoidal extensions built in Theorem \ref{rectext} has the least possible growth as $R \to \infty$.

\begin{theorem} \label{asym}
For every $R>1$, let $\Psi_R \in H^1(\Omega_{R})$ be a solution of \eqref{div}. Then
$$\liminf_{R\to\infty} \dfrac{\|\nabla \Psi_R\|_{L^2(\Omega_R)}+\|\Psi_R\|_{L^4(\Omega_R)}}{\sqrt[4]{R}} >0 \, .$$
\end{theorem}
\noindent
\begin{proof}
Let $\gamma > 0$ and consider the convex subset of $H^1_c(\Omega_{R})$ given by
$$
\mathcal{H}_{\gamma}(\Omega_{R}) = \{v\in H^1_c(\Omega_{R}) \ | \ v=\gamma \ \text{ on } \ \partial Q_{R}\}\, .
$$
Define
$$
I_\gamma(R) \doteq \inf\limits_{v \in \mathcal{H}_{\gamma}(\Omega_{R})} \left( \|\nabla v \|_{L^{2}(\Omega_{R})} + \|v\|_{L^{4}(\Omega_{R})} \right) \qquad \forall R > 1 \, ,
$$
and let $\varepsilon=\varepsilon(h) > 0$ be the unique positive solution of the algebraic equation
\begin{equation} \label{eps1}
\dfrac{2h}{3 \pi^{3}}\, \varepsilon^{4} + \dfrac{2\sqrt2}{h^{1/4}}\, \varepsilon =1\, .
\end{equation}
It suffices to prove that
\begin{equation} \label{eps2}
\liminf\limits_{R \to \infty} \dfrac{\, I_\gamma(R) \,}{\sqrt[4]{R}} \geq \varepsilon \gamma  \, .
\end{equation}

By homogeneity we can take $\gamma=1$, and \eqref{eps2} will be proved once we show that, for any possible choice of $V_{R} \in \mathcal{H}_1(\Omega_{R})$, it holds
	\begin{equation} \label{claimlimit}
	\liminf\limits_{R \to \infty} \dfrac{\|\nabla V_{R} \|^{4}_{L^{2}(\Omega_{R})} + \|V_{R}\|^{4}_{L^{4}(\Omega_{R})}}{R} \geq \varepsilon^{4} \, .
	\end{equation}
So, for any $R >1$, choose a function $V_{R} \in \mathcal{H}_1(\Omega_{R})$: two cases may occur. If
	$$
	\liminf\limits_{R \to \infty} \dfrac{\|\nabla V_R \|^{4}_{L^{2}(\Omega_{R})}}{R} \geq \varepsilon^{4} \, ,
	$$
	then the claimed limit \eqref{claimlimit} certainly holds. Otherwise, we have that
	$$
	\liminf \limits_{R \to \infty} \dfrac{\|\nabla V_R \|^{4}_{L^{2}(\Omega_{R})}}{R} = \delta \in \left[0, \varepsilon^{4} \right) \, .
	$$
	Then there exists a sequence of functions $V_{R_{k}} \in \mathcal{H}_1(\Omega_{R_k})$, $k \in \N$, where $(R_{k})_{k \in \mathbb{N}} \subset (1, \infty)$ is
	an increasing and divergent sequence, such that
	\begin{equation} \label{claimlimit3}
	\lim\limits_{k \to \infty} \dfrac{\|\nabla V_{R_k} \|^{4}_{L^{2}(\Omega_{R_k})}}{R_{k}} = \delta \, .
	\end{equation}
	By contradiction, assume that, possibly on a further subsequence $(R_j)_{j \in \mathbb{N}} \subset (R_k)_{k \in \mathbb{N}}$, it holds
	\begin{equation} \label{claimlimit2}
	\lim\limits_{j \to \infty} \dfrac{\|V_{R_{j}} \|^{4}_{L^{4}(\Omega_{R_j})}}{R_j} < \varepsilon^{4} - \delta \, .
	\end{equation}
	For every $j \in \mathbb{N}$ we then define a.e.\ in the rectangle $Q_j \doteq (-R_{j},R_{j}) \times (-h,h)$ the function
	$$
	W_{j}(x,y) =
	\left\lbrace
	\begin{array}{ll}
	1 - V_{R_{j}}(x,y) \quad & \text{if} \ \ (x,y) \in \Omega_{R_{j}} \, ,\\[5pt]
	1 \quad & \text{if} \ \ (x,y) \in \overline{B_{1}} \, ,
	\end{array}
	\right.
	$$
	so that $W_{j} \in H^{1}_{0}(Q_j)$. In view of \eqref{sobquad1} we then have, for all $j \in \N$,
	\begin{equation} \label{blowup1}
	\|\nabla V_{R_{j}} \|_{L^{2}(\Omega_{R_{j}})} ^4= \|\nabla W_{j} \|_{L^{2}(Q_j)}^4\geq \dfrac{3 \pi^{3}}{8} \dfrac{R^{2}_{j} + h^{2}}{R^{2}_{j} h^{2}} \| W_{j} \|^{4}_{L^{4}(Q_j)} = \dfrac{3 \pi^{3}}{8} \dfrac{R^{2}_{j} + h^{2}}{R^{2}_{j} h^{2}} \left[ \pi + \int_{\Omega_{R_j}} \! \! | 1 - V_{R_{j}} |^{4} \right].
	\end{equation}
	Then we notice that an application of H\"older's inequality yields
	\begin{equation} \label{holder2}
	\begin{aligned}
	\int_{\Omega_{R_j}} | 1 - V_{R_{j}} |^{4} & = |\Omega_{R_j} | - 4 \int_{\Omega_{R_j}} V_{R_{j}}  + \int_{\Omega_{R_j}} V_{R_{j}}^2 \left(6  - 4 V_{R_{j}} +  V_{R_{j}}^2 \right) \geq |\Omega_{R_j} | - 4 \|V_{R_{j}} \|_{L^{1}(\Omega_{R_j})} \\[6pt]
	& \geq |\Omega_{R_j} | - 4 |\Omega_{R_j} |^{3/4} \|V_{R_{j}} \|_{L^{4}(\Omega_{R_j})} = (4hR_{j} - \pi) - 4(4hR_{j} - \pi)^{3/4} \|V_{R_{j}} \|_{L^{4}(\Omega_{R_j})} \, .
	\end{aligned}
	\end{equation}
	Since we assumed \eqref{claimlimit2}, this shows that
	$$
	\begin{aligned}
	\liminf\limits_{j\to\infty} \dfrac{1}{R_{j}} \int_{\Omega_{R_j}} | 1 - V_{R_{j}} |^{4} & \geq 4 \left[h - (4h)^{3/4} \left( \varepsilon^{4} - \delta \right)^{1/4} \right] \geq 4 \left[ h - (4h)^{3/4} \varepsilon \right] .
	\end{aligned}
	$$
	Combining \eqref{eps1} with \eqref{blowup1}, this implies that
	$$
	\liminf \limits_{j \to \infty} \dfrac{\|\nabla V_{R_j} \|^{4}_{L^{2}(\Omega_{R_j})}}{R_j} \geq \dfrac{3 \pi^{3} }{2h^{2}} \left[ h - (4h)^{3/4} \varepsilon \right] = \varepsilon^{4} > \delta \, ,
	$$
	which contradicts \eqref{claimlimit3}.
\end{proof}

\begin{remark}
	The lower bound in \eqref{eps2} can be improved by taking a larger $\varepsilon$, both by noticing that \eqref{claimlimit} is merely a sufficient condition for \eqref{eps2} and by observing that in \eqref{holder2} one may write
	$$
	\int_{\Omega_{R_j}} V_{R_{j}}^2 \left(6  - 4 V_{R_{j}} +  V_{R_{j}}^2 \right) \geq \max \left\{ 2 \int_{\Omega_{R_j}} V_{R_{j}}^2 \ \  , \ \ \dfrac{1}{3} \int_{\Omega_{R_j}} V_{R_{j}}^4 \right\} \qquad \forall j \in \N \, .
	$$
	However, we will not pursue this further because we are mainly interested in the growth rate \eqref{eps2}.
\end{remark}

\begin{figure}[H]
	\begin{center}
		\includegraphics[height=58mm,width=100mm]{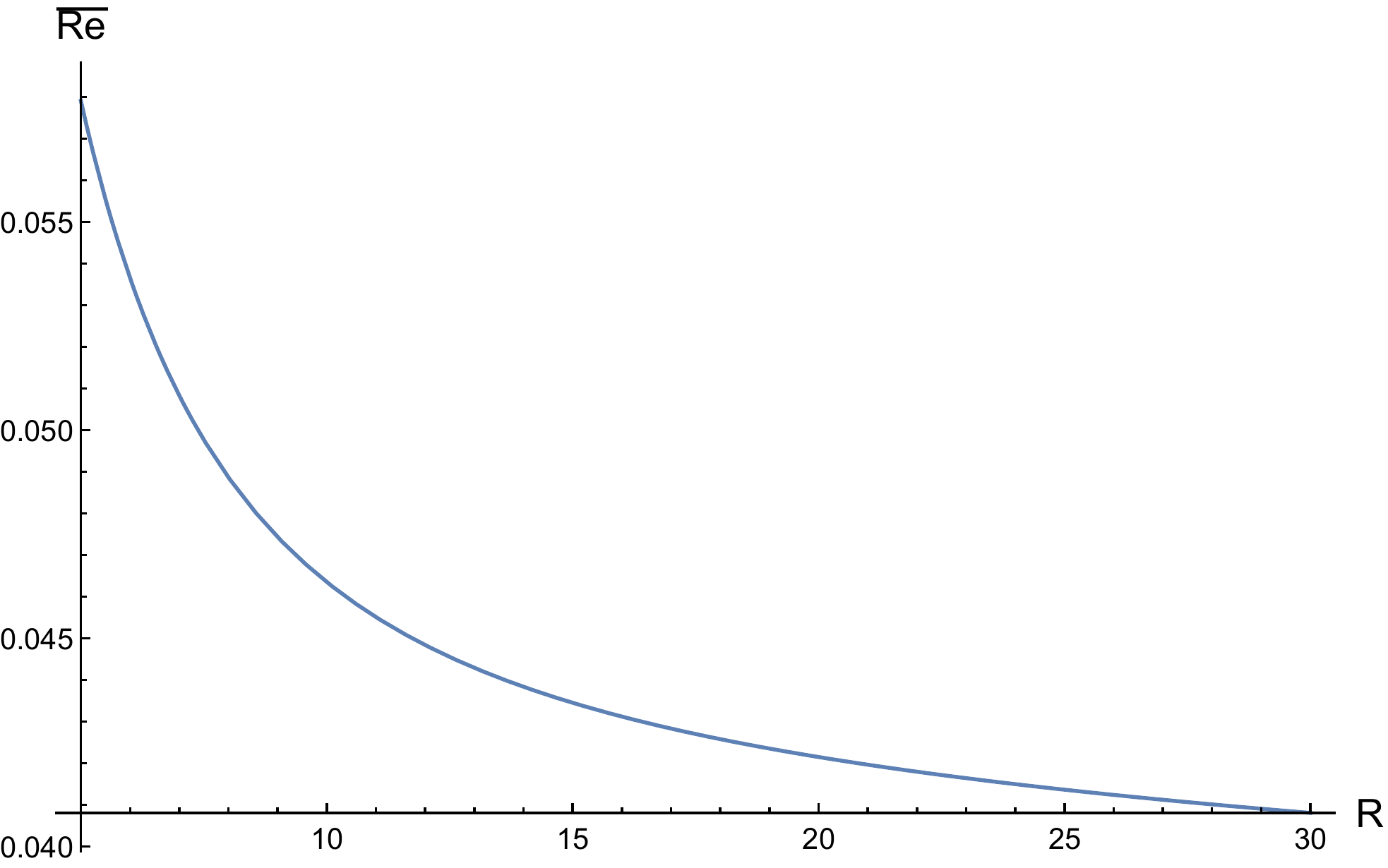}
	\end{center}
	\vspace*{-5mm}
	\caption{Upper bound for the Reynolds number \eqref{cotarey2} as a function of $R > h$, for $h=5$.}\label{reynolds2}
\end{figure}

If we maintain both the radius of the obstacle $B_1$ and the
height $2h$ of the outer rectangle fixed, the result of Corollary \ref{UVconstantcor} yields an explicit
upper bound for the Reynolds number that ensures the existence of a unique weak solution of problem \eqref{nsstokes0} which, moreover, displays the symmetry \eqref{simpro}. In Figure \ref{reynolds2} we plot the quantity
\begin{equation} \label{cotarey2}
\overline{\text{Re}} \doteq \dfrac{\dfrac{\pi}{2} \sqrt{\dfrac{3 \pi}{2}} \dfrac{\sqrt{R^2 + h^2}}{Rh} } {2 \sqrt{\mathcal{B}_{1} } + \sqrt{\dfrac{\pi}{2Rh}} \, \sqrt[4]{ \dfrac{3 \pi \mathcal{B}_{2}}{2}(R^2 + h^2)} } \, ,
\end{equation}
as a function of $R > h$, for some fixed value of $h > 1$.

\begin{remark} \label{psiopt}
By letting $R \to \infty$, we notice that \eqref{cotarey2} implies both
$$
\lim\limits_{h \to 1}  \overline{\text{Re}} = 0 \quad \forall R > 1 \qquad\mbox{and}\qquad
\overline{\text{Re}} \asymp \dfrac{1}{\sqrt[4]{R}} \quad \text{as} \quad R \to \infty \, .
$$
In fact, this asymptotic behavior as $R\to\infty$ cannot be improved. Indeed, the family of solenoidal extensions $\Psi_{R} \in H^{1}(\Omega_{R})$ built in Theorem \ref{rectext} satisfies
	$$
	\|\nabla\Psi_{R}\|_{L^{2}(\Omega_{R})} \ \ \text{does not depend on} \ R \qquad \text{and} \qquad \|\Psi_{R}\|_{L^{4}(\Omega_{R})} \sim \sqrt[4]{4hR} \quad \text{as} \quad R \to \infty \, .
	$$
	In view of Theorem \ref{asym} we deduce that the growth of these norms is ``minimal''.
\end{remark}

The results of the present section give a second connection between symmetry breaking in the stationary Navier-Stokes equations and Sobolev minimizers.
\begin{remark} \label{conn2}
	Theorem \ref{UVconstant0} states that, as long as uniqueness is guaranteed, symmetry of solutions of \eqref{nsstokes0} is observed. In turn, uniqueness is ensured provided that \eqref{umbral1}, which strongly depends on the Sobolev constant $\mathcal{S}_{R}$, holds. As we saw in Theorem \ref{tobias}, the minimizers for $\mathcal{S}_{R}$ lose symmetry as $R \to \infty$.
\end{remark}

\par\medskip\noindent
{\bf Acknowledgements.} The first Author is supported by the PRIN project {\em Direct and inverse problems for partial differential equations:
	theoretical aspects and applications} and by INdAM. The second Author is supported by the Research Programme PRIMUS/19/SCI/01, by the program GJ19-11707Y of the Czech National Grant Agency GA\v{C}R, and by the University Centre UNCE/SCI/023 of the Charles University in Prague.

\phantomsection
\addcontentsline{toc}{section}{References}
\bibliographystyle{abbrv}
\bibliography{references}
\vspace{5mm}

\begin{minipage}{80mm}
	Filippo Gazzola\\
	Dipartimento di Matematica\\
	Politecnico di Milano\\
	Piazza Leonardo da Vinci 32\\
	20133 Milan - Italy\\
	E-mail: filippo.gazzola@polimi.it
\end{minipage}
\hfill
\begin{minipage}{80mm}
	Gianmarco Sperone\\
	Department of Mathematical Analysis\\
	Charles University in Prague\\
	Sokolovská 83\\
	186 75 Prague - Czech Republic\\
	E-mail: sperone@karlin.mff.cuni.cz
\end{minipage}

\vspace{1cm}
\begin{minipage}{80mm}
	Tobias Weth\\
	Institut für Mathematik\\
	Goethe-Universität Frankfurt\\
	Robert-Mayer-Stra{\ss}e 10\\
	60629 Frankfurt - Germany\\
	E-mail: weth@math.uni-frankfurt.de
\end{minipage}
\end{document}